\documentclass[a4paper
]{amsart}

\usepackage[
left=40truemm,
right=40truemm
]{geometry}
\usepackage{tikz}
\numberwithin{figure}{section}

\usetikzlibrary{positioning,intersections, calc, arrows,decorations.markings, angles}
\colorlet{Green}{green!70!black!}
\usepackage{amsfonts, amssymb, amsmath, verbatim, amsthm, mathrsfs, amscd, enumerate, ascmac, fancyhdr, caption, framed, subfigure}

\usepackage{doi}

\usepackage{bbm}
\numberwithin{equation}{section}

\usepackage{hyperref}
\usepackage{bookmark}

\hypersetup{
	colorlinks=true,       
	linkcolor=blue,          
	citecolor=magenta,        
	filecolor=violet,      
	urlcolor=orange           
}

\theoremstyle{plain}
\newtheorem{theorem}{Theorem}[section]
\newtheorem{lemma}[theorem]{Lemma}
\newtheorem*{lemma*}{Lemma}
\newtheorem{proposition}[theorem]{Proposition}
\newtheorem*{proposition*}{Proposition}
\newtheorem{corollary}[theorem]{Corollary}
\theoremstyle{definition}

\newtheorem*{claim*}{Claim}

\newtheorem{conj}[theorem]{Conjecture}

\theoremstyle{remark}
\newtheorem*{remark}{Remark}

\def\sgn{\mathop{\mathrm{sgn}}\nolimits}

\def\d{\mathrm{d}}

\def\C{\mathcal{C}}

\def\1{\mathbbm 1}

\newcommand{\R}{{\mathbb R}}
\newcommand{\Z}{{\mathbb Z}}

\newcommand{\B}{{\mathbb B}}
\newcommand{\e}{\epsilon}
\def\j{\textbf{j}}

\def\proj{\textsf{Proj}}
\def\pst{{r_{\textsf{ST}}}}

\def\tang{\textsf{tang}}
\def\trans{\textsf{trans}}

\author[Kinoshita]{Shinya Kinoshita}
\address[Shinya Kinoshita]{Department of Mathematics, Institute of Science Tokyo, Meguro-ku, Tokyo, 152-8551, Japan}
\email{kinoshita@math.titech.ac.jp}

\author[Ko]{Hyerim Ko}
\address[Hyerim Ko]{Department of Mathematics, and Institute of Pure and Applied Mathematics, Jeonbuk National University, Jeonju, 54896, Republic of Korea}
\email{kohr@jbnu.ac.kr}

\author[Shiraki]{Shobu Shiraki}
\address[Shobu Shiraki]{Departamento de Matem\'atica, Instituto Superior T\'ecnico, Universidade de Lisboa, Av. Rovisco Pais 1, 1049-001 Lisbon, Portugal}
\email{shobushiraki@tecnico.ulisboa.pt}

\thanks{
}

\begin{document}
\date{\today}
\title[
]
{
Maximal estimates for orthonormal systems of wave equations
}

\subjclass[2010]{}
\begin{abstract}
This paper investigates maximal estimates of the wave operators for orthonormal families of initial data. We extend the classical maximal estimates for the wave operator by making partial progress on maximal estimates for orthonormal systems in low dimensions. Our novel approach is based on a geometric analysis of the kernel of wave operators within the framework of Schatten $2$ estimates. In particular, we exploit Wolff's geometric lemma on the intersection patterns of thickened spheres.
\end{abstract}

\maketitle


\section{Introduction}
\subsection{Motivation}
Let $n\ge1$. For $s>0$, let $H^s(\R^n)$ denote the inhomogeneous Sobolev space of order $s$ equipped with the norm $\|g\|_{H^s}=\| (1-\Delta)^{s/2}g\|_{L^2}$. For $g\in H^s(\R^n)$, we consider the wave equation $\partial_t^2u-\Delta u=0$ with initial conditions $u(x,0)=g(x)$ and $\partial_tu(x,0)=0$.
The solution $u$, denoted by $\cos(t\sqrt{-\Delta})g$, can be represented as
\[
\cos (t \sqrt{-\Delta}) g(x) = \frac 12 \big(U_t+U_{-t}\big)g(x),
\]
where $U_tg(x)$ is the solution to the half-wave equation
\begin{equation}\label{wave-equation}
i \partial_t u + \sqrt{- \Delta} u=0  
\end{equation}
under the initial condition $u(x,0)=g(x)$.
Furthermore, $U_tg(x)$ has the explicit representation
\[
U_tg(x)
=
(2\pi)^{-n}\int_{\R^n}e^{i(x\cdot \xi+t|\xi|)} \widehat g(\xi)\,\d\xi.
\]

Maximal estimates for various operators have drawn considerable attention from different perspectives. A classical result due to \cite{Cowling, Walther} establishes that maximal estimates for single wave operator:
\begin{equation}\label{e:max single wave}
    \Big\|\sup_{0<t<1} |U_t g|\Big\|_{L^2(\B^n)}
    \lesssim
    \|g\|_{H^s(\R^n)}
\end{equation}
holds if and only if $s>\frac12$. Here, $\mathbb{B}^n(x,r)$ denotes the $n$-dimensional ball of radius $r>0$ centered at $x$ and $\mathbb{B}^n=\mathbb{B}^n(0,1)$ represents the unit ball.

Historically, local maximal estimates of this type were first explored in the context of the Schr\"odinger equation, famously known as Carleson's problem \cite{Carleson,DK82}. While our focus here is on the wave equation, interested readers may refer to \cite{Bourgain16, DGL17, DZ19} for recent breakthroughs on the Schr\"odinger equation variant.
The analysis of \eqref{e:max single wave} for the wave equation is notably more straightforward, following from relatively simple arguments, whereas the Schr\"odinger case requires more delicate recent techniques.
For a detailed discussion, see Appendix~\ref{appendix-single}.

The maximal estimates for orthonormal systems $(f_j)_j$ in $H^s(\R^n)$ present a natural extension of \eqref{e:max single wave}, where we consider infinitely many initial data, rather than a single datum. This generalization can be formulated as 
\begin{equation}\label{main}
\Big\|\sup_{0<t<1} \sum_j \lambda_j|U_t f_j|^2 \Big\|_{L^{\beta}(\mathbb{B}^n)}
\lesssim
\|\lambda\|_{\ell^\beta}.
\end{equation}
Here, $(\lambda_j)_j\in \ell^\beta(\mathbb{C})$ and $\beta\geq1$.
In fact, this is equivalent to
\[
\Big\|  
\sup_{0<t<1} \big( \sum_j |U_t g_j|^2 \big)^{\frac 12} \Big\|_{L^{2\beta} (\mathbb{B}^n)}
\lesssim \Big( \sum_j \|g_j\|_{H^s}^{2\beta} \Big)^{\frac{1}{2\beta}}
\]
for orthogonal systems $(g_j)_j$ in $H^s(\R^n)$, and it is straightforward to see \eqref{main} holds for $\beta=1$ if $s>\frac12$ as a consequence of the aforementioned classical result for a single datum, along with a simple application of the triangle inequality. Furthermore, if the regularity is sufficiently smooth, namely when $s>\frac n2$, the Bessel inequality guarantees \eqref{main} holds with $\beta=\infty$, hence for any $\beta\geq1$ by the embedding $\ell^\beta \subset \ell^\infty$. In fact, assuming $(g_j)_j$ is orthonormal in $L^2(\R^n)$, and thus so is $(\widehat{g}_j)_j$, we observe  
\begin{align}\label{e:application Bessel}
    \sum_j \frac{|\lambda_j|}{\|\lambda_j\|_{\ell^\infty}}|U_t\langle  D\rangle^{-s} g_j(x)|^2
    \leq
    \sum_j \Big|\int_{\R^n} e^{i(x\cdot\xi+t|\xi|)}\langle \xi\rangle^{-s} \widehat{g_j}(\xi)\,\d\xi\Big|^2
    \leq
    \|\langle \cdot \rangle^{-s}\|_{L^2}^2,
\end{align}
where $\langle D\rangle:=(1-\Delta)^{1/2}$ and $\langle \xi \rangle = (1+|\xi|^2)^{1/2}$. The right-hand side remains bounded as long as $s>\frac n2$.  
Consequently, this also shows that there is nothing particularly interesting in the case $n=1$. 
Additionally, when $s=\frac n2$, the estimate \eqref{main} holds for $\beta\in[1,\infty)$, as a consequence of \eqref{tri-interpol} below together with a minor adjustment at low frequency.

With those discussions in mind, in the current manuscript, we aim to establish \eqref{main} for the largest possible $\beta$ (depending on $s$) for $n\geq2$ and $s\in(\frac12,\frac n2)$ by making use of the orthogonal structure of the data.

\subsection{Previous works}
The study of Strichartz estimates for orthonormal system has been made by Frank, Lewin, Lieb, and Serringer \cite{FLLS} and Frank--Sabin \cite{FS_AJM, FS_Survey} for some operators such as Schr\"odinger operators. 
The result was later extended by Bez et al. \cite{BHLNS} up to
a sharp range for non-sharp admissible cases considering homogeneous initial datum.
Strichartz estimates for various dispersive operators (cf. \cite{BLN_Forum})
and some of the boundary cases have been obtained by \cite{BKS_TLMS, BKS_RIMS}. 

The Strichartz estimates for orthonormal system for wave equation seem more difficult and the sharp characterization in terms of the function spaces remains an open problem. In particular, the following diagonal Strichartz estimate has been conjectured (for a more general statement, see \cite{BLN_Forum, BKS_RIMS}).
\begin{conj}\label{c:Strichartz}
Let $n\geq2$ and $q\geq \frac{2(n+1)}{n-1}$. Then, we have
\[
\Big\| \sum_j \lambda_j |U_t f_j|^2\Big\|_{L_{x,t}^\frac q2(\R^{n+1})}
\lesssim
\|\lambda\|_{\ell^\beta}
\]
for orthonormal systems $(f_j)_j$ in the homogeneous Sobolev space $\dot{H}^s(\R^n)$ if $\beta\leq \frac{n}{2(n+1)}q$.
\end{conj}

The study of maximal estimates in the context of orthonormal systems is even more recent. Bez, Lee, and Nakamura \cite{BLN_Selecta} were the first to obtain a maximal-in-time estimate for the Schr\"odinger equation in one dimension, achieving the largest possible $\beta$ when the regularity is at its minimal value, $s = \frac{1}{4}$. Their method involves first establishing corresponding maximal-in-space estimates and then applying the technique of swapping the roles of space and time to derive the maximal-in-time estimates.
Subsequently, Bez and two of the present authors \cite{BKS_TLMS} adopted a more direct approach to extend these results beyond $s = \frac{1}{4}$, including for fractional Schr\"odinger operators. They also addressed cases where $s$ is relatively close to $\frac{n}{2}$ in higher dimensions.

\subsection{Main results}
Similar challenges to those in the Strichartz estimates arise when considering the maximal estimates for wave operators. We address these challenges and obtain new, non-trivial results that represent partial progress toward the optimal bounds.

To begin, we propose our conjecture concerning the maximal estimates for orthonormal systems in the wave equation framework. The conjecture is based on modifications of the counterexamples from the single-datum setting in \cite{CHL, BKS_RIMS}, following a similar approach to that used in \cite{BKS_TLMS} for the Schr\"odinger equation. 
Throughout the paper, we consider the orthonormal function of initial data $(f_j)_j$ in $H^s(\mathbb{R}^n)$. 
However, all our results remain valid for $(f_j)_j$ in the homogeneous Sobolev space $\dot{H}^s(\R^n)$ by a similar computation to \eqref{e:application Bessel}.
\begin{conj}\label{conj-wave}
Let $n\geq2$ and $\frac12< s < \frac n2$. 
The estimate \eqref{main} holds for all orthonormal system $(f_j)_j$ in $H^s(\R^n)$ if and only if
\[
\beta\le \min\Big\{ \frac n{n-2s}, \frac {n-1}{n+1-4s} \Big\}.
\]
\end{conj}
One might attempt to approach the problem by adapting the arguments used to establish single-datum results, such as \eqref{e:max single wave} 
based on the Sobolev embedding in order to address Conjecture~\ref{conj-wave}. However, this would depend on Conjecture~\ref{c:Strichartz}, and even if the full Conjecture~\ref{c:Strichartz} were proven,  it would still not be sufficient. We will elaborate on this point in Section~\ref{sec:loc-reduct}.
We prove \eqref{main} with a nontrivial range of $\beta=\beta(s)$ when $n=2,3,4$.

\begin{theorem}\label{t:d34}
Let $n=3,4$ and $\frac12 < s< \frac n2$. Then \eqref{main} holds for all orthonormal system $(f_j)_j\in H^s(\R^n)$ provided
\begin{align}\label{betas1}
\beta< \min \Big \{ \frac{2n-1}{2(n-2s)}, \frac{2n-3}{2n-1-4s}\Big\}.
\end{align}
\end{theorem}

See Figure~\ref{f:3dim} below to visualize the improvement compared to Conjecture~\ref{conj-wave} when $n=3$. The argument employed in Theorem~\ref{t:d34} extended to the case $n=2$, yielding the bound
$\beta<\min\{ \frac{5}{8-4s}, \frac{3}{7-4s}\}$
(see the Subsection~\ref{n2small} for further discussions).

Moreover, by exploiting the self-similar structure of the inequality through a certain conic decomposition, we achieve a further improvement beyond the $n=2$ case of Theorem~\ref{t:d34}, as illustrated in Figure~\ref{f:2dim}.

\begin{theorem} \label{t:refiment d2}
Let $n=2$ and $\frac12<s<1$. Then, \eqref{main} holds for all orthonormal system $(f_j)_j\in H^s(\mathbb R^2)$ provided
\begin{align}\label{betas2}
\beta<\min\left\{\frac{9}{14(1-s)},\frac{5}{12-14s}\right\}.
\end{align}
\end{theorem}

\medskip

\textbf{Novelty of the paper.}
The main obstacle to proving Theorem~\ref{t:d34} and \ref{t:refiment d2} lies in establishing the sharp estimate \eqref{main} for $\beta=2$. After a series of reductions, the problem reduces to proving a rather simple bilinear estimate with the optimal decay rate $a>0$, namely, 
\begin{align}\label{e:bilinear a}
\iint\iint h(x,t) h(x',t') \chi_{\mathbb{V}_{k}}(x-x',t-t')\,\d x\d t\d x'\d t'
\lesssim
2^{-ak}\|h\|_{L_x^2L_t^1}^2
\end{align}
for all $h \in L_x^2L_t^1$ where $\mathbb V_k$ denotes the one-sided cone $\{(\xi,|\xi|)\in \R^n\times \R\}$ of thickness approximately $2^{-k}$.

\begin{figure}[htbp]
\centering
\begin{tikzpicture}[scale=7]
	\draw [->] (1/2,1.03) node (yaxis) [above] {\small $\frac 1\beta$}|- (3/2+0.05,0) node ( xaxis) [right] {$s$}  ;
\draw [->] (1/2,0)--(1/2,1.03);
\draw (1/2,0)--(1/2,1)--(3/2,1)--(3/2,0)--(1/2,0);
\draw [densely dotted] (1/2,1/2)--(3/4,1/2)--(3/4,0);
\draw [densely dotted] (3/4,1/2)--(7/8,1/2)--(7/8,0);

\path [name path=conj1] (3/2,0)--(3/4,1/2);
\path [name path=conj2] (1,0)--(1/2,1);
\path [name path=theorem1] (5/4,0)--(1/2,1);
\path [name path=theorem2] (3/2,0)--(1/4,1);
\path [name path=theorem3] (2,0)--(1/4,1);

\path [name path=conj5] (3/2,0)--(.05,1/2);

\path[name intersections={of=conj1 and conj2, by={A}}];
\path[name intersections={of=theorem1 and theorem2, by={B}}];

\fill [teal, opacity=0.2] (1/2,0)--(3/2,0)--(A)--(1/2,1)--(1/2,0);
\fill [violet, opacity=0.2] (3/2,0)--(3/2,1)--(1/2,1)--(B)--(3/2,0);

\draw [teal, thick] (3/2,0)--(3/4,1/2)--(1/2,1);
\draw [violet, thick] (3/2,0)--(7/8,1/2)--(1/2,1);
\draw [gray, densely dashed, thick] (3/2,0)--(1/2,1);
    
\node [below] at (1/2,0) {\scriptsize $\frac12$};
\node [left] at (1/2,0) {\scriptsize $0$};
\node [below] at (7/8,0) {\scriptsize $\frac78$};
\node [below] at (3/4,0) {\scriptsize $\frac 34$};
\node [below] at (6/4,0) {\scriptsize $\frac 32$};
\node [left] at (1/2,1) {\scriptsize $1$};
\node [left] at (1/2,1/2) {\scriptsize $\frac12$};
	
\end{tikzpicture}
\caption{The relationship between $s$ and $\frac{1}{\beta}$ when $n = 3$. The green region indicates the range of parameters for which inequality~\eqref{main} is conjectured to fail, as formulated in Conjecture~\ref{conj-wave}. In contrast, the purple region corresponds to the range where inequality~\eqref{main} is proved to hold, as established in Theorem~\ref{t:d34}.}
\label{f:3dim}

\centering
\begin{tikzpicture}[scale=7]

\draw (1/2,0)--(1,0)--(1,1)--(1/2,1)--(1/2,0);
\draw [->] (1/2,0)--(1/2,1.03);
\draw [->] (1/2,0)--(1.05,0);
\node [right] at (1.05,0) {$s$};
\node [above] at (1/2,1.025) {$\frac1\beta$};

\path [teal, name path=conj1] (1,0)--(0,1);
\path [teal, name path=conj2] (3/4,0)--(1/2,1);

\path [violet, name path=theorem1] (6/7,0)--(1/2,1);
\path [violet, name path=theorem2] (1,0)--(5/14,1);
\path [name path=theorem3] (1.5,0)--(1/4,1);

\path[name intersections={of=conj1 and conj2, by={A}}];

\path[name intersections={of=theorem1 and theorem2, by={B}}];

\fill [teal, opacity=0.2] (1/2,0)--(1,0)--(A)--(1/2,1)--(1/2,1)--(1/2,0);
\fill [violet, opacity=0.2] (1,0)--(1,1)--(1/2,1)--(B)--(1,0);

\draw [densely dotted] (1/2,1/3) node [left] {\scriptsize $\frac 13$}--(2/3,1/3)--(2/3,0) node [below] {\scriptsize $\frac 23$};

\node [below] at (1/2,0) {\scriptsize $\frac12$};
\node [left] at (1/2,0) {\scriptsize $0$};
\node [below] at (1,0) {\scriptsize $1$};
\node [below] at (1/2,0) {\scriptsize $\frac12$};
\node [left] at (1/2,1) {\scriptsize $1$};

\draw [densely dotted] (B)--(B|-,1);
\draw [densely dotted] (B)--(1,1/2);
\node [right] at (1,1/2) {\scriptsize $\frac12$};
\node [above] at (B|-,1) {\scriptsize $\frac{19}{28}$};

\draw [teal, thick] (1,0)--(2/3,1/3)--(1/2,1);
\draw [violet, thick] (1,0)--(19/28,1/2)--(1/2,1);
\draw [gray, densely dashed, thick] (1,0)--(1/2,1);

\end{tikzpicture}
\caption{The relationship between $s$ and $\frac{1}{\beta}$ when $n = 2$. The green region indicates the range of parameters for which inequality~\eqref{main} is conjectured to fail, as formulated in Conjecture~\ref{conj-wave}. In contrast, the purple region corresponds to the range where inequality~\eqref{main} is proved to hold, as established in Theorem~\ref{t:refiment d2}.}
\label{f:2dim}
\end{figure}

 It is worth noting that a similar estimate is also crucial for establishing Strichartz estimates for the wave operators. However, applying the same strategy used for the Schr\"odinger operator (as in \cite{BLN_Forum, BKS_TLMS}) only yields trivial estimates for \eqref{e:bilinear a}. (We discuss this in detail in Section~\ref{sec:kernel}.)
The novelty of this paper lies in our investigation of nontrivial estimates for \eqref{e:bilinear a}, achieved by considering the intersection of the spheres in $\R^{n+1}$.


\subsection{Pointwise convergence}
Naturally, the equation \eqref{wave-equation} can be extended to the operator-valued equation: 
\begin{equation}\label{eq:HalfWaveOp}
i \partial_t \gamma = [-\sqrt{- \Delta} , \gamma]
\end{equation}
under the condition $\gamma(0)=\gamma_0$,
where $\gamma_0$ is a compact self-adjoint operator on $L^2(\R^n)$. 
Precisely, taking $\gamma_0$ as the projection onto the span of $f \in L^2(\R^n)$, the solution $\gamma(t)$ of \eqref{eq:HalfWaveOp}, described as $\gamma(t) = U_t \gamma_0 U_{-t}$, is the projection onto the span of $U_t f$. 

An immediate consequence of Theorems~\ref{t:d34} and \ref{t:refiment d2} is the pointwise convergence of the density of an operator in a certain class. For $\beta\geq1$ let 
 $\C^{\beta}$ be the Schatten space of all compact operators on $L^2(\mathbb R^{d+1})$ equipped with the norm $\|T\|_{\C^{\beta}}:=\|\lambda\|_{\ell^{\beta}}<\infty$. Here, $(\lambda_j)_j$ are the singular values of $T$. 
Let us recall a precise definition of a density of an operator in the Schatten--Sobolev space $\mathcal C^{\beta,s}$ defined by
\[
\mathcal{C}^{\beta,s} = \{ \gamma \in \mathrm{Com}(H^{-s}(\R^n), H^s(\R^n)) \, : \, \| \gamma\|_{\mathcal{C}^{\beta, s}} :=\| \langle -\Delta \rangle^{\frac{s}{2}} \gamma \langle -\Delta \rangle^{\frac{s}{2}} \|_{\mathcal{C}^{\beta}(L^2(\R^n))} < \infty \}.
\]
Here, $\mathrm{Com}(H^{-s}(\R^n), H^s(\R^n))$ denotes the set of compact operators from $H^{-s}(\R^n)$ to $H^{s}(\R^n)$.

For $\beta,s$ satisfying the conditions in \eqref{betas1} or \eqref{betas2}, the density of an operator in $\mathcal{C}^{\beta,s}$ can be represented as a limit of a sequence of density of finite-rank operators. 
Precisely, let $\gamma_0$ be an operator in $\mathcal{C}^{\beta,s}$ where $\beta,s$ satisfy the conditions in \eqref{betas1} or \eqref{betas2}. Suppose that $T:= \langle -\Delta \rangle^{\frac{s}{2}} \gamma \langle -\Delta \rangle^{\frac{s}{2}}$ is a self-adjoint operator in $\mathcal{C}^{\beta}(L^2(\R^n))$. Then there exists an orthonormal family $(f_j)_j \subset L^2(\R^n)$ and $(\lambda_j)_j \in \ell^{\beta}(\R_+)$ where $\R_+$ is the set of positive real numbers such that 
\[
T^N := \sum_{j=1}^N \lambda_j \Pi_{f_j} \to T \qquad ( N \to \infty) 
\]
in $\mathcal{C}^{\beta}(L^2(\R^n))$. 
Here $\Pi_{f}$ denotes the orthogonal projection onto the span of $f$. The density of $T^N$ is given by 
\[
\rho_{T^N}(x) =  \sum_{j=1}^N \lambda_j |f_j(x)|^2,
\]
and then the density of $\gamma^N := \langle -\Delta \rangle^{-\frac{s}{2}} T^N \langle -\Delta \rangle^{-\frac{s}{2}}$ is formally given by 
\[
\rho_{\gamma^N}(x) = \sum_{j=1}^N \lambda_j |\langle -\Delta \rangle^{-\frac{s}{2}} f_j(x)|^2.
\]
We define the density of $\gamma$ as a limit of $\rho_{\gamma^N}$, namely,
\[
\rho_{\gamma}(x) = \sum_{j} \lambda_j |\langle -\Delta \rangle^{-\frac{s}{2}} f_j(x)|^2.
\]
Since the solution of \eqref{eq:HalfWaveOp} with initial data $\gamma_0$ can be written as $\gamma(t)=U_t \gamma_0 U_{-t}$, the density of $\gamma(t)$ takes the form
\[
\rho_{\gamma(t)} (x)= \rho_{\raisebox{-2pt}{\scriptsize{$U_t \gamma_0 U_{-t}$}}}(x) = \sum_{j} \lambda_j |\langle -\Delta \rangle^{-\frac{s}{2}}U_t f_j(x)|^2.
\]
For further details, we refer to Section 6 of \cite{BKS_TLMS}.

Consider an analogue of classical pointwise convergence problems with finding the largest possible class of initial data $\gamma_0$ for which the following holds:
\begin{align}\label{pt-wise}
\lim_{t\rightarrow 0} \rho_{\gamma(t)}(x)=\rho_{\gamma_0}(x)
\end{align}
for almost every $x \in \R^n$. In view of Conjecture~\ref{conj-wave}, we expect that \eqref{pt-wise} holds provided
\[
s >
\max\Big\{ \frac{n}{2} - \frac{n}{2\beta},~
\frac{n+1}{4} - \frac{n-1}{4\beta}
\Big\}.
\]
As a consequence of Theorem~\ref{t:d34} and \ref{t:refiment d2}, we obtain the following.

\begin{corollary}\label{t:PWC}
Let $s\in (\frac12,\frac n2)$ and $1\le \beta \le \infty$. Then we have
\begin{enumerate}
\item [$($i\,$)$]
For $n=3,4$, if $\gamma_0 \in \mathcal{C}^{\beta,s}$ is self-adjoint with $\beta, s$ satisfying 
\[
s>\max\big\{ \frac n2-\frac{2n-1}{4\beta}, ~\frac {2n-1}4-\frac{2n-3}{4\beta}\big\},
\]
then $\rho_{\gamma(t)}$ and $\rho_{\gamma_0}$ are well-defined and \eqref{pt-wise} holds.
\vspace{.1cm}
\item [$($ii\,$)$] For $n=2$, if $\gamma_0 \in \mathcal{C}^{\beta,s}$ is self-adjoint with $\beta$, $s$ satisfying 
\[
s>\max\big\{1-\frac{9}{14\beta},~ \frac67-\frac{5}{14\beta}\big\},
\]
then $\rho_{\gamma(t)}$ and $\rho_{\gamma_0}$ are well-defined and \eqref{pt-wise} holds.
\end{enumerate}
\end{corollary}

Following the approach in \cite{BLN_Selecta}, Corollary~\ref{t:PWC} follows from Theorem~\ref{t:d34}.
For a detailed proof, we refer to \cite{BLN_Selecta, BKS_TLMS}.


\subsection*{Organization of the paper}
In Section~\ref{s:necessary}, we discuss the necessary conditions for the Conjecture~\ref{conj-wave}. 
Section~\ref{sec:loc-reduct} reduces Theorem~\ref{t:d34} and \ref{t:refiment d2} to Proposition~\ref{p:frec loc beta2}, which addresses the case where $\beta=2$. 
In Section~\ref{sec:further reduct}, we further reduce Proposition~\ref{p:frec loc beta2} to Proposition~\ref{prop:integral}, which involves geometric estimates such as \eqref{e:bilinear a}. 

In Section~\ref{s:d34}, we prove Proposition~\ref{prop:integral} for dimensions $n=3,4$ employing Wolff's geometric estimates for sphere intersection (Lemma~\ref{l:wolff} and \ref{l:wolff higher dim}). For dimension $n=2$, we provide a refined proof in Section~\ref{s:refinement}, exploiting the self-similar structure of the main estimates.

Appendix~\ref{appendix-single} reviews known results on the maximal estimates for the wave and Klein--Gordon operators for single initial datum, providing detailed proofs of almost sharp estimates up to endpoints. In Appendix~\ref{Appendix:B}, we present a refined proof of Wolff's lemma (Lemma~\ref{l:wolff} and \ref{l:wolff higher dim}) for the reader's convenience.

\medskip

\subsection*{Notational conventions}

Throughout the paper, we use the following notations.

\begin{itemize}

\item $\varphi \in C_0^\infty(\R)$: a smooth function that is compactly supported on $(\frac12,2)$ such that $\varphi_k=\varphi(2^{-k}\cdot)$, $k\in\Z$ forms a partition of unity.
\item $P_k$: the standard Littlewood--Paley projection operator associated with a smooth function $\varphi_k$.

\item $\chi_A$: the characteristic function of the set $A$.

\item $A \lessapprox B$:  $A\le C 2^{\epsilon k} B$ for given $\epsilon>0$, $k\ge1$, real numbers $A,B>0$, and some constant $C>0$.

\item $\mathbb{B}^n(x,r)$: the $n$-dimensional ball of radius $r>0$ centered at $x$. Also, $\mathbb{B}^n=\mathbb{B}^n(0,1)$.

\item  $\mathbb{V}_{k,l}(x_0,t_0)$: a $2^{-k}$-neighborhood of an upper cone centered at $(x_0,t_0)$ with height $2^{-l}$, namely, $\{(x,t): ||x-x_0|-(t-t_0)|\leq 2^{-k},\ 0\leq t-t_0\leq 2^{-l}\}$. We also write
$\mathbb{V}_{k,l}=\mathbb{V}_{k,l}(0,0)$.

\item $\mathbb{W}(w_1,w_2)$: $\proj_{\R^n}( \mathbb{V}_{k,l}(w_1)\cap \mathbb{V}_{k,l}(w_2))$ for given points $w_1,w_2\in \R^{n+1}$.

\end{itemize}

\section{Necessary conditions}\label{s:necessary}

In this section, we discuss the necessary conditions for $s,\beta$ for which \eqref{main} holds.
\begin{proposition}\label{prop:necessary}
Let $n\ge2$ and $\beta \ge1$. Then there exists orthonormal families of functions $\{f_j\}$ in $H^s(\R^n)$ for which the estimate \eqref{main} fails when
\[
\beta
> 
\min\Big\{ \frac n{n-2s}, \frac {n-1}{n+1-4s} \Big\}.
\]
\end{proposition}

\subsection{Proof of $\beta\le \frac n{n-2s}$} 
Suppose \eqref{main} holds. To show the necessity of $\beta \le \frac{n}{n-2s}$, 
we employ the same example as in \cite{BKS_TLMS}. For the sake of completeness, we include a sketch of the proof. For an integer $N\gg1$, let 
\[
J_j
=
[N^{-1}j,N^{-1}j+(100N)^{-1}]
\]
for $j\in \mathbb Z$.
One may find $\eta\in C_0^\infty(\mathbb R^n)$ supported on the unit ball $\mathbb B^n$ such that $\widehat{\eta}$ is non-negative, and $\int_{|\xi|\geq1}\widehat{\eta} (\xi)\,\d \xi\leq \e \|\widehat \eta\|_1$ for a small $\e>0$. 
Let $\mathit{J}:=\{\mathbf j=(j_1,\dots,j_n)\in \mathbb Z^n: J_{j_i} \cap \mathbb B^1 \neq \emptyset~\text{for all~} i=1,\dots,n\}$. Then, for 
$\j \in \mathit{J}$,
set the orthonormal initial data
\[
f_{\j}(x)
=
c N^{\frac n2-s} \eta \big(10N(x-N^{-1}\j)\big)
\]
for some appropriate constant $c>0$ so that $\|f_{\mathbf j}\|_{H^s}=1$. 
We also choose $\lambda_{\j}=1$ if $\j\in \mathit{J}$ and $0$ otherwise so that $\|\lambda\|_{\ell^\beta}\sim |\mathit{J}|^\frac{1}{\beta}\sim N^\frac{n}{\beta}$. 

With our choice of $(f_{\textbf j})_{\textbf j}$ and $(\lambda_{\j})_{\j}$, one may observe that for each $\j$,
\begin{align*}
|U_t f_{\j}(x)| 
=
cN^{\frac n2-s} \Big| \int e^{i(10N(x-N^{-1}\j)\cdot \xi+10 Nt|\xi|)} \widehat{\eta}(\xi)\,\d \xi\Big|
\gtrsim 
N^{\frac n2-s}
\end{align*}
whenever $t=t_N:=(100N)^{-1}$ and $x\in \mathbb{B}^n(N^{-1}\j, (100N)^{-1}) =: \B_\j^n$
since the phase is fairly small (at most $\frac12$).
Therefore, we obtain
\[
\Big\| \sum_{\j\in\mathit{J}} |U_t f_{\j}|^2\Big\|_{L_x^1L_t^\infty}
\geq 
\sum_{\j\in \mathit{J}} \|U_{t_N}f_{\mathbf j}\|_{L_x^2(\B_\j^n)}^2
\gtrsim
N^{n-2s}.
\]
Hence, assuming \eqref{main}, we obtain $N^{n-2s}\lesssim N^\frac{n}{\beta}$, which only holds for large $N$ if $\beta \le \frac{n}{n-2s}$.
\qed

\subsection{Proof of $\beta\le \frac{n-1}{n+1-4s}$}

Next, we shall show $\beta\le \frac{n-1}{n+1-4s}$ 
is necessary for \eqref{main} to hold. For an integer $N\gg 1$, we set
\[
\bar{\mathit{J}_j}
:=
[N^{-\frac12} j ,N^{-\frac12} j+100^{-1}N^{-\frac12}]
\]
for $j \in \mathbb Z$.
Let $\phi_P$ be a Schwartz function essentially supported on 
\[
P
=
\{\xi=(\xi_1,\bar{\xi})\in\mathbb R\times \mathbb R^{n-1} : N\leq \xi_1\leq 2N,\ |\bar{\xi}|\leq N^{\frac12}\}.
\]
For $\bar{\j}\in \overline{\mathit{J}}$ where $\overline{\mathit J}:=\{\bar{\j} = (j_1,\dots,j_{n-1}) \in \mathbb Z^{n-1}: \bar{J_{j_i}}\cap \mathbb B^1\not=\emptyset ~\text{for all}~ i=1,\dots,n\}$, we define an orthonormal initial data
\[
\widehat{f_{\,\bar{\j}}}(\xi)
=
c N^{-\frac{n+1}{4}-s} e^{-i \langle \bar{\j},\, \bar{\xi}\rangle} \phi_P(\xi),
\]
where the constant $c>0$ is chosen so that $\|f_{\,\bar{\mathbf j}}\|_{H^s}=1$ and $(f_{\,\bar{\mathbf j}})_{\bar{\mathbf j}}$ forms an orthonormal basis in $H^s$. We also choose $\lambda_{\,\bar{\j}}=1$ if $\bar{\j}\in\overline{\mathit{J}}$ and $0$ otherwise so that $\|\lambda_{\,\bar {\mathbf j}}\|_{\ell^\beta}=|\overline{\mathit{J}}|^\frac1\beta\lesssim N^{\frac{n-1}{2\beta}}$.
With our choice of $(f_{\,\bar{\j}})_{\bar{\j}}$ and $(\lambda_{\,\bar{\j}})_{\bar{\j}}$, one may observe that for each $\bar{\j}$ 
\[
|U_tf_{\,\bar{\j}}(x,t)|
=
cN^{-\frac{n+1}{4}-s}\Big|\int e^{i(\langle \bar{x}-\bar{\j},\,\bar{\xi}\rangle+(x_1+t)\xi_1+t(|\xi|-\xi_1))}\phi_P(\xi)\,\d\xi\Big|
\gtrsim
cN^{\frac{n+1}{4}-s}
\]
whenever $t=-x_1$ and $x\in Q_{\bar{\j}}:=\{x\in \B^n(0,100^{-1}) : |\bar{x}-\bar{\j}|\leq 100^{-1} N^{-\frac12}\}$ since $|\xi|-\xi_1=|\bar{\xi}|^2/(|\xi|+\xi_1)\leq |\bar{\xi}|^2/\xi_1$ guarantees the phase is small enough. Therefore, noting $Q_{\bar{\j}}$ are disjoint and each has the volume $\sim N^{-\frac{n-1}{2}}$, we obtain
\[
\Big\| \sum_{\bar{\j}} |U_t f_{\,\bar{\j}}|^2\Big\|_{L_x^1L_t^\infty}
\gtrsim
\sum_{\bar{\mathbf j}} \|U_{-x_1}f_{\,\bar{\mathbf j}}\|_{L_x^2(Q_{\bar{\j}})}^2
\gtrsim
N^{\frac{n+1}{2}-2s} \Big|\bigcup_{\bar{\j}} Q_{\,\bar{\j}}\Big|
\sim
N^{\frac{n+1}{2}-2s}.
\]
Hence, assuming \eqref{main}, we obtain $N^{\frac{n+1}{2}-2s}\lesssim N^{\frac{n-1}{2\beta}}$, which only holds for large $N$ if $\beta \le \frac{n-1}{n+1-4s}$. 
\qed

\section{Localization and reductions}\label{sec:loc-reduct}

In this section, we introduce several tools to demonstrate Theorem~\ref{t:d34} and \ref{t:refiment d2}, some of which are fairly standard but included for the sake of completeness.

\subsection{Frequency localization}
We start with a standard argument, also employed in \cite[Section 6.2]{BKS_TLMS}. In the rest of the paper, we denote by $(f_j)_j$ an orthonormal basis in $H^s(\mathbb{R}^n)$ and by $(\lambda_j)_j$ the sequence in $\ell^\beta$ for appropriate $\beta\geq1$ in the context. 
Let $\varphi\in C_0^\infty(\R)$ supported in $\{r\in\mathbb R: 2^{-1}\leq r<2\}$ and $\varphi_k=\varphi(2^{-k}\cdot)$ for $k\in \Z$. Associated with $\varphi$ we define 
the standard Littlewood--Paley projection operator $P_k$ by $P_kf=(\varphi_k(|\cdot|)\widehat{f})^\vee$.

Then, Theorem~\ref{t:d34} and \ref{t:refiment d2} follow from the frequency localized estimates
\begin{equation}\label{main frequency loc}
\Big\|\sup_{0<t<1} \sum_j \lambda_j|U_tP_k f_j|^2 \Big\|_{L_x^{\beta}(\mathbb{B}^n)}
\lesssim
2^{2\varsigma k}
\|\lambda\|_{\ell^\beta}
\end{equation}
for $k\geq1$,
where $\varsigma>\varsigma^\star$ given by
\[
\varsigma^\star =
\begin{cases}
\max\big\{1-\frac{9}{14\beta},~ \frac67-\frac{5}{14\beta}\big\}, &\quad n=2,\\[1ex]
\max \big\{\frac n2-\frac{2n-1}{4\beta},~ \frac {2n-1}{4}-\frac {2n-3}{4\beta} \big\}, &\quad n=3,4.
\end{cases}
\]

In fact, for $s\in(0,\frac n2)$ the Littlewood--Paley decomposition yields
\begin{align*}
    \Big\|
    \sup_{0<t<1} \sum_j \lambda_j|U_t\langle D\rangle^{-s} f_j|^2 
    \Big\|_{L_x^{\beta}(\mathbb{B}^n)}
    \leq 
    \sum_{k\in \Z}
    2^{-2sk}
    \Big\|
    \sup_{0<t<1} \sum_j \lambda_j|U_t P_k f_j|^2 
    \Big\|_{L_x^{\beta}(\mathbb{B}^n)}.
\end{align*}
The term with $k<0$ can be controlled using the trivial bound
\begin{equation}\label{tri-infty}
\Big\|\sum_j \lambda_j|U_t P_k f_j|^2\Big\|_{L_x^\infty L_t^\infty}
\lesssim
2^{nk}
\|\lambda\|_{\ell^\infty}
\end{equation}
which is obtained by Bessel's inequality (refer to \eqref{e:application Bessel} or \cite{FS_AJM, BHLNS}).

Because of this, we restrict our attention to $k \geq 0$ in the rest of the paper.

\subsection{Interpolation and reduction to $\beta=2$}\label{sec:reduct} 
Showing \eqref{main frequency loc} when $\beta=1$ and $\beta=\infty$ is relatively easy.
For \eqref{main frequency loc} with $\beta=1$, we recall the classical result for a single initial datum (see \eqref{e:max single wave} with $p=2$ and Appendix~\ref{appendix-single}). This combined with triangle inequality gives
\begin{equation}\label{tri-1}
\Big\|\sum_j \lambda_j|U_t P_k  f_j|^2\Big\|_{L_x^1L_t^\infty}
\lesssim
\sum_j \lambda_j \big\|U_t P_k  f_j\|_{L_x^2L_t^\infty}^2
\lesssim
2^k
\|\lambda\|_{\ell^1}.
\end{equation}
By interpolating between \eqref{tri-1} and \eqref{tri-infty}, we obtain
\begin{align}\label{tri-interpol}
\Big\|\sup_{0<t<1} \sum_j \lambda_j|U_t\langle D\rangle^{-s} P_k f_j|^2 \Big\|_{L_x^{\beta}(\mathbb{B}^n)}
\lesssim
2^{(n-\frac{n-1}{\beta}-2s)k}
\|\lambda\|_{\ell^\beta}
\end{align}
for all $1\le \beta\le \infty$. Summing over $k\geq0$, the estimate \eqref{main} holds whenever $s>\frac n2-\frac{n-1}{2\beta}$, but this is significantly larger than the conjectured regularity $\max\{\frac n2-\frac{n}{2\beta},\frac{n+1}{4}-\frac{n-1}{4\beta}\}$.

A typical approach to the problem is to establish some (ideally, sharp) bounds for $\beta=2$, \eqref{est-c2}
and interpolate the cases $\beta=1,2,\infty$. 
Due to Proposition~\ref{prop:necessary}, the best one can hope for is that 
\begin{align}\label{est-c2}
\Big\| \sup_{0<t<1} \sum_j \lambda_j |U_tP_kf_j|^2 \Big\|_{L_x^2(\mathbb{B}^n)} 
\lessapprox
2^{\frac{1}{2}(n+1-\sigma)k} \|\lambda\|_{\ell^2}
\end{align}
for $\sigma=\sigma_n$ given by 
\begin{align*}
\sigma_n:=
\begin{cases}
	\frac12, \quad \text{if} \quad n=2, \\
	1, \quad \text{if} \quad n\ge3.
\end{cases}
\end{align*}
Here, 
we write $A\lessapprox B$ to mean $A \lesssim C_\epsilon^{\epsilon k} B$ for $A,B>0$, and some constant $C_\epsilon>0$ depending on arbitrary small $\epsilon>0$. Although the bound \eqref{est-c2} alone is not sufficient to obtain sharp estimates for all $\beta$ via interpolation, this becomes relevant only after \eqref{est-c2} is established with the optimal exponent $\sigma = \sigma_n$. Further refinement is then required, as observed for instance in \cite{BKS_TLMS} in the context of the Schr\"odinger equation.

The trivial bound \eqref{tri-interpol} coincides with \eqref{est-c2} when $\sigma=0$.
As we mentioned in the Introduction, a natural approach would be to adapt the techniques from the single-datum setting, such as \eqref{e:max single wave} or \eqref{maximal single modified}, to handle this conjecture. However, this approach proves inefficient in this context, as even resolving Conjecture~\ref{c:Strichartz} would still only give us \eqref{est-c2} with $\sigma = \frac{n}{2(n+1)}$. Therefore, in the remainder of this paper, we will instead employ a geometric argument to obtain the following results. To state our results, set \begin{equation}\label{estar}
\sigma_n^\star
=
\begin{cases}
\frac27=\frac{1}{4}+\frac{1}{28}\ &\text{if} \quad n=2,\\
\frac12  \ &\text{if} \quad n=3,4.
\end{cases}
\end{equation}

\begin{proposition}\label{p:frec loc beta2}
Let $n=2,3,4$. Then for any $\epsilon>0$, \eqref{est-c2} holds for $\sigma=\sigma_n^\star$.
\end{proposition}

Consequently, our goal \eqref{main frequency loc} is obtained by interpolating this with \eqref{tri-1} and \eqref{tri-infty}. To obtain Proposition~\ref{p:frec loc beta2}, we make use of some geometric properties of cones and circles (Section~\ref{s:d34}). Moreover, exploiting the self-similar structure of the inequality, we further refine our results for $n=2$.

\medskip

\subsection{Duality}
Recall that $\C^{p}$ is the Schatten space
equipped with the norm $\|T\|_{\C^{p}}:=\|\lambda\|_{\ell^{p}}<\infty$ where $(\lambda_j)_j$ are the singular values of $T$. 
The most important property of the Schatten space in the current manuscript is that $\|\cdot\|_{\C^2}$ is the Hilbert--Schmidt norm and coincides with the norm $\|\cdot\|_{L^2(\mathbb R^{d+1}\times\mathbb R^{d+1})}$ of the kernel of the corresponding operator. 

We apply the following lemma to dualize \eqref{main frequency loc} to the corresponding bound with respect to $\C^2$-norm.
\begin{lemma}[Duality principle \cite{FS_AJM, BHLNS}]\label{l:duality}
Suppose $T$ is a bounded linear operator from $L^2(\mathbb{R}^n)$ to $L_x^{r}L_t^{q}$ for some $q,r \geq2$. Also, let $\beta\geq1$.  Then, 
\[
\Big\|\sum_j\lambda_j|Tf_j|^2 \Big\|_{L_x^{\frac {r}{2}}L_t^{\frac {q}{2}}}
\lesssim
\|\lambda\|_{\ell^{\beta}}
\]
holds for all orthonormal systems $(f_j)_j$ in $L^2(\mathbb R^n)$ and all sequences $(\lambda_j)_j \in \ell^{\beta}(\mathbb C)$ if and only if 
\[
\|WTT^*\overline{W}\|_{\C^{\beta'}}
\lesssim 
\|W\|_{L_x^{\widetilde{r}}L_t^{\widetilde{q}}}^2
\]
for all $W\in L_x^{\widetilde{r}}L_t^{\widetilde{q}}$. Here, $\widetilde{\cdot}$ denotes the half conjugate given by
\[
\frac{1}{p}+\frac{1}{\widetilde{p}}=\frac12
\]
and $\cdot'$ denotes the (usual H\"older) conjugate given by
\begin{align*}
\frac1\beta+\frac{1}{\beta'}=1.
\end{align*}
\end{lemma}

By Lemma~\ref{l:duality}, Proposition~\ref{p:frec loc beta2} is equivalent to the following: Let $T_k=(U_tP_k)(U_tP_k)^*$, then
\begin{align}\label{dual-c2}
\|WT_k\overline{W}\|_{\C^{2}}
\lesssim
2^{\frac{(n+1-\sigma)}2k}
\|W\|_{L^{4}_xL_t^2}^2
\end{align}
holds for all $W\in L_x^{4} L_t^2$ with $\sigma=\sigma_n^{\star}$ when $n=2,3,4$.
Here,
$T_kF=F\ast K_k$
where the kernel $K_k$ is given by
\begin{align}\label{kernel-k}
K_k(x,t)
=
\int e^{i(x\cdot \xi+t|\xi|)}\varphi_k^2(|\xi|)\,\d \xi.
\end{align}

As we mentioned at the beginning of this section, $\mathcal C^2$ is the Hilbert--Schmidt class and the $\mathcal C^2$-norm can be interpreted as
\begin{equation}
\begin{aligned}\label{C2-integral}
&\|WT_k\overline W\|_{\mathcal C^2(L^2(\mathbb B^{n+1}),L^2(\mathbb B^{n+1}))}^2\\
&=
\iint_{\B^{n+1}}\iint_{\B^{n+1}}
\Big| W(x,t) K_k(x-x',t-t') W(x',t')\Big|^2\,\d x \d t \d x' \d t'.
\end{aligned}
\end{equation}
Here, we may assume that $W(x,t)$, $\overline W(x,t)$ are supported on $\mathbb B^{n+1}$ since $U_t P_kf$ is supported on $\mathbb B^{n}\times (0,1)$.

\subsection{Kernel estimates}\label{sec:kernel}
We recall the well-known property of the kernel of the wave operator (see for example \cite{Lee03}), which can be derived by using polar coordinates and asymptotic behavior of the Bessel functions.

\begin{lemma}\label{l:dispersive}
Let $k\ge1$ and $(x,t)\in \R^n\times \R$. Then for any $N\ge1$, 
\[
\Big|
\int e^{i(x\cdot \xi+t|\xi|)} \varphi_k^2(|\xi|)\,\d \xi
\Big|
\le C_N
2^{nk} (1+2^k|x|)^{-\frac{n-1}{2}} (1+2^k||x|-|t||)^{-N}.
\]
\end{lemma}

If one invokes Lemma~\ref{l:dispersive} in a trivial manner, it follows from \eqref{C2-integral} that
\[
\|WT_k\overline W\|_{\mathcal C^2}^2
\lesssim 2^{2nk}\int |W(x,t)|^2|\overline W(x',t')|^2(1+2^k|x-x'|)^{-(n-1)}\,\d x\d x'\d t\d t'.
\]
A simple application of Young's convolution inequality, for instance, immediately yields \eqref{dual-c2} with 
 \[
 \sigma=0, \quad n \ge2.
 \]
However, this bound can be trivially obtained by an interpolation between the estimates \eqref{tri-1} and \eqref{tri-infty}.

Unlike the Schr\"odinger operator, the argument presented above does not yield sharp $\mathcal C^2$ estimates.
It appears that the use of the additional properties of the kernel $K_k$, which is supported on a small neighborhood of the cone $\{(x,|x|): x\in \mathbb{B}^n \}$, is crucial for proving the sharp results. 
To overcome this difficulty, we exploit the support property of the kernel to derive the improved $\mathcal C^2$-estimates which are to be stated in the subsequent sections.

\section{Further reduction to geometric estimates}\label{sec:further reduct}
As briefly mentioned in the introduction, a bilinear integral formed by a geometry from the wave equation plays a central role in understanding Proposition~\ref{p:frec loc beta2}.  
For $0\le l \le k$, let 
\[
\mathcal I_{k,l} 
:=
\iint h(x,t) h(x',t')\chi_{\mathbb{V}_{k,l}}(x-x',t-t')\,\d x \d x' \d t \d t',
\]
where 
$\mathbb{V}_{k,l}$ denotes an upper cone in $\R^{n+1}$ defined by
\[
\mathbb{V}_{k,l}
:=
\big\{(x,t)\in \mathbb B^{n}\times [0,1]: ||x|-|t||\lesssim 2^{-k},~ 2^{-l-1}\le |x|\le 2^{-l} \big\}.
\]
Then, we claim the following.
\begin{proposition}\label{prop:integral}
Let $n=3,4$ and $0\le l  < k$. For $\sigma=\sigma_n^\star$ given by \eqref{estar}, then we have
\begin{align}\label{integral2}
\mathcal I_{k,l}
\lessapprox 
2^{-\sigma k}2^{-(n-1)l} \|h\|_{L_x^2L_t^1}^2
\end{align}
holds for any $h \in L_x^2L_t^1$.

\end{proposition}

We observe below that 
Proposition~\ref{p:frec loc beta2} follows from Proposition~\ref{prop:integral}. 
The proof of Proposition~\ref{prop:integral} will be postponed to Section~\ref{s:d34}.

\subsection{Reduction from Proposition~\ref{p:frec loc beta2} to Proposition~\ref{prop:integral}}
By the reduction in \eqref{dual-c2}, it suffices to show that 
\begin{align}\label{e:C2}
\|WT_{k}\overline{W}\|_{\C^2}^2
&\lessapprox
2^{(n+1-\sigma_n^\star)k} \|W\|_{L_x^4L_t^2}^4
\end{align}
where $\sigma_n^\star$ is given by \eqref{estar}.

Recalling the kernel $K_k$ defined by \eqref{kernel-k},
we further apply a spatial decomposition, given by $(2^l)_l$, to the kernel, namely, for a fixed $k$, 
we define $\chi_l(t)=\chi_{[2^{-l-1},2^{-l})}(t)$ for $0 \le l \le k-1$ and $\chi_l(t)=\chi_{[0,2^{-k})}(t)$ for $l=k$.
We set $K_{k,l}(x,t)=\chi_l(|x|)K_k(x,t)$
 so that 
\[
T_kF=\sum_{0 \le l \le k}T_{k,l}F:=\sum_{0 \le l \le k}F \ast K_{k,l}.
\]
Let us denote $h(x,t):=|W(x,t)|^2$ so that 
\begin{align}\label{WTKLW}
\|WT_{k,l}\overline{W}\|_{\C^2}^2
&=\!\!
\iiiint h(x,t) h(x',t') |K_{k,l}(x-x',t-t')|^2\,\d x \d x'\d t\d t'.
\end{align}

To analyze the integral, we consider the support of the kernel.
When $l=k$, Lemma~\ref{l:dispersive} implies that the kernel $|K_{k,l}|$ is bounded by $2^{nk}\chi_{\mathbb B^{n+1}(0,2^{-k})}$, which is easy to handle.
In fact, by setting $B_k:=\mathbb B^{n+1}(0,2^{-k})$, Young's convolution inequality yields
\[\|h (h\ast \chi_{B_k})\|_{L_{x,t}^1} 
\lesssim \|h\|_{L_x^2L_t^1}\|h\ast \chi_{B_k}\|_{L_x^2L_t^\infty}
\lesssim 
2^{-nk}
\|h\|_{L_x^2L_t^1}^2.
\] 
This estimate is sufficient to prove Proposition~\ref{prop:integral} with the conjectured value $\sigma=\sigma_n=1$ in the case $l=k$.

When $0\le l\le k-1$, we decompose the domain into three parts: two corresponding to 
\[
\pm \mathbb{V}_{k,l}^\epsilon
:=
\big\{(x,\pm t)\in \mathbb B^{n}\times [0,1]: ||x|-|t||\lesssim 2^{-(1-\epsilon)k},~ 2^{-l-1}\le |x|\le 2^{-l} \big\},
\]
which denote the upper and lower cones, 
and their complement 
$\mathbb B^{n+1}\setminus \big\{\mathbb{V}_{k,l}^\epsilon \cup (-\mathbb{V}_{k,l}^\epsilon) \big\}$. The complement is a minor part since, by Lemma~\ref{l:dispersive}, the kernel satisfies
$$
|K_{k,l}(x,t)\chi_{\mathbb B^{n+1}\setminus \{ \mathbb{V}_{k,l}^\epsilon \cup (-\mathbb{V}_{k,l}^{\epsilon})\}}(x,t)|\lesssim_N 2^{-Nk}
$$
for any $N\ge1$.

This estimate suffices to consider the integral \eqref{WTKLW} over $\mathbb{V}_{k,\ell}^\epsilon \cup (-\mathbb{V}_{k,l}^\epsilon)$.
We set
\[
\mathcal I_{k,l}^{\epsilon,\pm}
=
\iint h(x,t) h\ast \chi_{(\pm \mathbb{V}_{k,l}^\epsilon)}(x,t)\,\d x\d t.
\]
Since $|K_{k,l}(x,t)|^2\lesssim 2^{(n+1)k}2^{-(n-1)l}$ provided $(x,t) \in \pm \mathbb V_{k,l}^\epsilon$ by Lemma~\ref{l:dispersive}, combining three cases yields
\begin{align}\label{c2-Ikl}
\| WT_{k,l}\overline W\|_{\mathcal C^2}^2
\lesssim 
2^{(n+1)k} 2^{(n-1)l} \sum_\pm \mathcal I_{k,l}^{\epsilon,\pm} + 2^{-Nk}\|h\|_{L_x^2L_t^1}^2.
\end{align}

By symmetry, it is enough to estimate $\mathcal I_{k,l}^{\epsilon,+}$.
Further allowing $2^{\epsilon k}$ loss, we reduce our goal to estimating $\mathcal I_{k,l}$, which is the statement of Proposition~\ref{prop:integral}.

\section{Proof of Proposition~\ref{prop:integral} when $n=3,4$}\label{s:d34}

In this section, we complete the proof of Theorem~\ref{t:d34} establishing Proposition~\ref{prop:integral}. The key step involves analyzing the interaction between two $2^{-k}$-fattened cones in $\R^{n+1}$. In particular, the size of the vertically projected intersection plays a crucial role.

\subsection{Reduction to geometric intersection estimates}
By the H\"older inequality, we observe that
\begin{align}\label{Ikl-WW}
\mathcal I_{k,l} \le \| h \|_{L_x^2L_t^1}
\big\| h \ast \chi_{\mathbb{V}_{k,l}} \big\|_{L_x^2L_t^\infty}.
\end{align}
For simplicity, we denote by $w=(x,t)\in \R^n\times \R$.
After taking a square and expanding the integral of the right-hand side, and then applying Fubini's theorem, we see that
\begin{align}\label{e:sup inside}
    \|h*\chi_{\mathbb{V}_{k,l}}\|_{L_x^2L_t^\infty}^2
    &\leq
    \iint h(w_1) h(w_2) \int \sup_t\chi_{\mathbb{V}_{k,l}(w_1)\cap \mathbb{V}_{k,l}(w_2)}(x,t)\,\d x \d w_1 \d w_2 \nonumber\\
    &\leq
    \iint h(w_1) h(w_2) 
    |\mathbb{W}(w_1,w_2)|
    \, \d w_1 \d w_2,
\end{align}
where $|\mathbb{W}(w_1,w_2)|$ is the Lebesgue measure of the set defined by 
\begin{align}\label{Wn=proj}
\mathbb{W}(w_1,w_2) = \proj_{\R^n}( \mathbb{V}_{k,l}(w_1)\cap \mathbb{V}_{k,l}(w_2))
\end{align}
and $\mathbb{V}_{k,l}(w):=\mathbb{V}_{k,l}+w$, the translation of $\mathbb{V}_{k,l}$ by the vector $w$.
By dividing the intersection into level set $\{t=t_n\}$, it is now important to observe the intersection of the circles and spheres of thickness $2^{-k}$.


\subsection{Geometric estimates for sphere intersection}
We now measure the size of the intersection of two thickened circles.
For $x\in \R^n$ and $r>0$, we denote a $\delta$-neighborhood of the sphere $r\mathbb S^{n-1}$ centered at $x$ in $\R^n$ by
\[
\mathbb O_\delta(x,r)=\big\{ z\in \R^n: ||z-x|-r|\le \delta\big\}.
\]

The following is a slight variant of Wolff's lemma \cite[Lemma 3.1]{Wolff} (see also \cite{Marstrand} and \cite[Lemma B.2]{CDK}) which includes both the external and internal tangent cases. Define 
$\varDelta=\varDelta(x_1,x_2,r_1,r_2)$ by
\begin{align}\label{d2}
\varDelta
:=
\big||x_1-x_2|-|r_1+ r_2|\big| \times
\big||x_1-x_2|-|r_1-r_2| \big|.
\end{align}
\begin{lemma}\label{l:wolff}
Let $x_1,x_2 \in \R^2$ and $r_1,r_2>0$, and
$0<\delta \le 1$. 
Then,
\begin{equation}\label{CC refined}
    |\mathbb{O}_\delta(x_1,r_1)\cap \mathbb{O}_\delta(x_2,r_2)|
    \lesssim
    \frac{\delta^{\frac32} (r_1+r_2)}{(|x_1-x_2|+\delta)^{\frac12}}\left(\frac{\delta (r_1+r_2)}{\varDelta+\delta (r_1+r_2)}\right)^\frac12.
\end{equation}
\end{lemma}

The proof of Lemma~\ref{l:wolff} is almost identical to that of \cite{Wolff}, but stated in a more precise manner. For the reader's convenience, we provide the detail in Appendix~\ref{Appendix:B}.

There is also a natural extension of this result to higher dimensions, specifically for estimating the intersection of two spheres in $\R^n$ with $n\geq3$. This was observed in \cite[Lemma B.1]{CDK} by employing an additional slicing argument alongside Lemma~\ref{l:wolff}.
For its proof, see Appendix~\ref{Appendix:B} as well.

\begin{lemma}\label{l:wolff higher dim}
Let $n\ge3$ and $x_1,x_2\in \R^n$.
Let $r_1,r_2>0$. Then 
\begin{align}\label{intersect high}
    |\mathbb{O}_\delta(x_1,r_1)\cap \mathbb{O}_\delta(x_2,r_2)|
    \lesssim 
    \frac{\delta^2 (r_1+r_2)^{\frac32}}{|x_1-x_2|+\delta}
    \Big(\frac{\varDelta+\delta(r_1+r_2)}{|x_1-x_2|+\delta}\Big)^{\frac{n-3}2}.
\end{align}
\end{lemma}

\subsection{Proof of Proposition~\ref{prop:integral} when $n=3,4$}
We claim that 
\begin{align}\label{e:two cones3}
\sup_{0<t_1,t_2<1}
|\mathbb{W}(x_1,t_1,x_2,t_2)|
\lesssim 
2^{-k}2^{-(n-\frac 12)l} |x_1-x_2|^{-\frac{n-1}2}
\end{align}
provided $|x_1-x_2|\ge 2^{-k}$ when $n\ge3$. 
First, assuming \eqref{e:two cones3} for the moment, we proceed to prove \eqref{e:C2}. Note that $|x_1-x_2|\lesssim 2^{-l}$ by the definition of $\mathbb{V}_{k,l}$.
By H\"older's inequality, the right-hand side of \eqref{e:sup inside} is bounded by
\begin{align*}
2^{-k}2^{-(n-\frac 12)l} \| h\|_{L_x^2L_t^1} 
\big\| \|h(\cdot,t)\|_{L_t^1} \ast \big(|\cdot|^{-\frac{n-1}2}\chi_{|\cdot|\le 2^{-l}}\big)(x) \big\|_{L_{x}^2}.
\end{align*}
  By Young's convolution inequality, it follows that
\begin{align*}
\big\|h \ast \chi_{\mathbb{V}_{k,l}} \big\|_{L_x^2L_t^\infty}
\lesssim 
	2^{-\frac{k}{2}}2^{-\frac {3}{4}nl} \|h\|_{L_x^2L_t^1}.
\end{align*}

Combining this with \eqref{c2-Ikl}, the bound \eqref{Ikl-WW} yields \eqref{e:C2} with $\sigma_{n}^\star=\frac12$ when $n=3,4$ since $2^{-\frac {3n}4l} \le 2^{-(n-1)l}$ for $n\le 4$. This completes the proof of Proposition~\ref{prop:integral}.\\

It now remains to prove \eqref{e:two cones3}.
Let $0\le l <k$. 
We slice the space by the level sets of height $2^{-k}$, defined as $\{2^{-k} m\le t < 2^{-k}(m+1)\}$.
By \eqref{intersect high}, for any $1\le m \lesssim 2^{k-l}$ and $n\ge3$, and since $\varDelta+\delta(r_1+r_2)\lesssim 2^{-2l}$, we obtain that
\[
\big|\proj_{\R^n}\big(\mathbb{V}_{k,l}(w_1)\cap \mathbb{V}_{k,l}(w_2)\cap \{2^{-k} m\le t < 2^{-k}(m+1)\}\big)\big|
\]
is bounded by
\[
2^{-2k}2^{-(n-\frac32)l} |x_1-x_2|^{-\frac{n-1}2}.
\]
Summing over all $m$ gives \eqref{e:two cones3}.\qed

\subsection{Analogous results for $n=2$}\label{n2small}
A similar approach to that used for Proposition~\ref{prop:integral} in the case $n=3,4$ above also works in $n=2$.
In fact, one can obtain \eqref{integral2} with $\sigma=\frac14$.

By using Lemma~\ref{l:wolff}, instead of Lemma~\ref{l:wolff higher dim}, we have
\begin{align*}
\sup_{0<t_1,t_2<1}|\mathbb{W}(x_1,t_1,x_2,t_2)|
\lesssim
2^{-\frac k2} 2^{-2l}|x_1-x_2|^{-1/2}
\end{align*}
instead of \eqref{e:two cones3}.
By applying H\"older's inequality and Young's convolution inequality, we obtain
\[
\big\| h \ast \chi_{\mathbb{V}_{k,l}} \big\|_{L_x^2L_t^\infty}^2 
\lesssim 
2^{-\frac k2}2^{-\frac 72l}\|h\|_{L_x^2L_t^1}^2.
\]

This together with \eqref{Ikl-WW} gives \eqref{integral2} with $\sigma=\frac14$ since $2^{-\frac 74l} \le 2^{-l}$.
\qed

\smallskip

\begin{remark}
    When Lemma~\ref{l:wolff} was applied above, it was used in the form
    \[
    |\mathbb{O}_\delta(x_1,r_1)\cap \mathbb{O}_\delta(x_2,r_2)|
    \lesssim
    \frac{\delta^{\frac32} (r_1+r_2)}{(|x_1-x_2|+\delta)^{\frac12}}.
    \]
    In the forthcoming section, we apply Lemma~\ref{l:wolff} in a more precise manner and observe \eqref{integral2} holds with the order $\sigma=\sigma_2^\star=\frac14+\frac{1}{28}$, slightly improved over $\frac14$. This improvement is achieved not only through a more detailed analysis of the relevant geometry, but also by modifying the estimate to take advantage of a certain bootstrapping argument. This is only possible in $n=2$, as the volume of the intersection of two cones typically becomes large when they are tangent in this dimension. This point will be clarified further in the discussion below.
\end{remark}

\section{Refinement in two dimensions}\label{s:refinement}

The purpose of this section is to prove Theorem~\ref{t:refiment d2}, based on the reductions carried out earlier. To this end, we establish the following proposition. As before, we occasionally write $w = (x,t) \in \mathbb{R}^n \times \mathbb{R}$.

\begin{proposition}
Let $n=2$ and $0\le l < k$. Then,
\begin{align}\label{goal1}
\mathcal I_{k,l} 
\lessapprox
2^{-\sigma_2^\star k} 2^{-\frac{12}7l}\|h\|_{L_x^2L_t^1}^2
\end{align}
holds for any $h \in L_x^2L_t^1$, where $\sigma_2^\star=\frac14+\frac{1}{28}$ is given as in \eqref{estar}.
\end{proposition}
From \eqref{Ikl-WW}, it suffices to show that
\begin{align}\label{goal2}
\big\| h \ast \chi_{\mathbb{V}_{k,l}} \big\|_{L_x^2L_t^\infty}
\lessapprox
2^{-\sigma_2^\star k}
2^{-(\frac 74-\frac1{28})l} \|h\|_{L_x^2L_t^1}.
\end{align}

\subsection{Proof of the main Proposition}
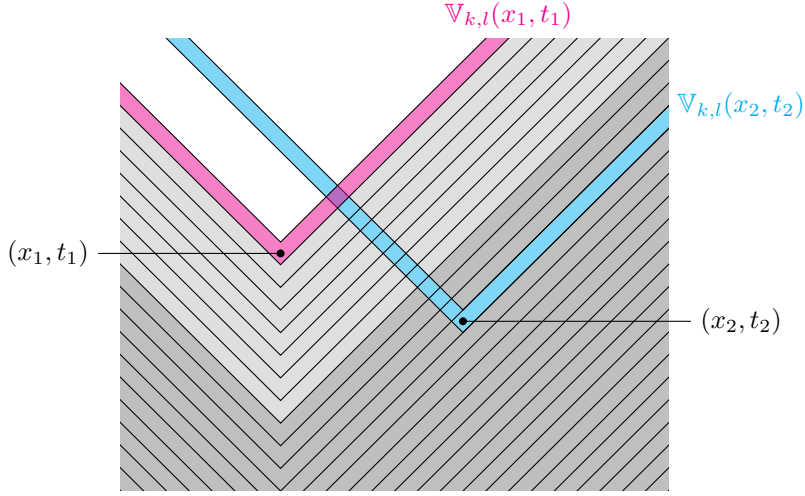
\begin{figure}[htbp]

\begin{tikzpicture}[scale=3]

\begin{scope}
  \clip (-0.7,-1) -- (-0.7,1) -- (1.7,1) -- (1.7,-1) -- cycle;

  \foreach \y in {1,2,...,7}{
    \begin{scope}[shift={(0,0.0 - 0.1*\y)}]
      \fill[lightgray,opacity=0.5] (2,2)--(0,0)--(-1,1)--(-1,1.1)--(0,0.1)--(2,2.1);
    \end{scope}
  }
    
  \foreach \y in {8,9,...,30}{
    \begin{scope}[shift={(0,0.0 - 0.1*\y)}]
      \fill[gray,opacity=0.5] (2,2)--(0,0)--(-1,1)--(-1,1.1)--(0,0.1)--(2,2.1);
    \end{scope}
  }

\begin{scope}[shift={(0.8,-0.3)}]

  \fill[white,opacity=1] (2,2)--(0,0)--(-2,2)--(-2,2.1)--(0,0.1)--(2,2.1);
  \fill[cyan,opacity=0.5] (2,2)--(0,0)--(-2,2)--(-2,2.1)--(0,0.1)--(2,2.1);
  \draw (2,2)--(0,0)--(-2,2); 
  \draw (2,2+0.1)--(0,0+0.1)--(-2,2+0.1);
\end{scope}

  \foreach \y in {1,2,...,30}{
    \begin{scope}[shift={(0,0.0 - 0.1*\y)}]
      \draw (2,2)--(0,0)--(-1,1); 
    \end{scope}
  }

  \fill[magenta,opacity=0.5] (1.6,1.6)--(0,0)--(-1,1)--(-1,1.1)--(0,0.1)--(1.6,1.7);
  \draw (1.6,1.6)--(0,0)--(-1,1); 
  \draw (1.6,1.6+0.1)--(0,0+0.1)--(-1,1+0.1);

\end{scope}

\fill [black] (0,0+0.05) circle [radius=1/2pt];
\draw (0,0+0.05)--(-0.8,0.05);
\node [left] at (-0.8,0.05) {$(x_1,t_1)$};
\fill [black] (4/5,-0.3+0.05) circle [radius=1/2pt];
\draw (4/5,-0.3+0.05)--(1.8,-0.3+0.05);
\node [right] at (1.8,-0.3+0.05) {$(x_2,t_2)$};

\node [magenta,above] at (1,1) {$\mathbb{V}_{k,l}(x_1,t_1)$}; 
\node [cyan,right] at (1+0.7,0.7) {$\mathbb{V}_{k,l}(x_2,t_2)$}; 

\end{tikzpicture}

\caption{An illustration of space-time decomposition. For fixed $w_1=(x_1,t_1)$, decompose the space-time into conical regions translated in time by $2^{-k}j$. We say that the cones close to $\mathbb{V}_{k,l}(x_1,t_1)$ are \textit{tangential} (light gray and $\mathbb V_{k,l}(x_1,t_1)$) and the cones away from $\mathbb{V}_{k,l}(x_1,t_1)$ are \textit{transversal} (gray). We distinguish the cases depending on whether the apex $(x_2,t_2)$ lies in a tangential or transversal cone.}
\label{f:decomposition}

\end{figure}

Here, we establish \eqref{goal2}. We begin by recalling \eqref{e:sup inside}, which is
\begin{align*}
\|h*\chi_{\mathbb{V}_{k,l}}\|_{L_x^2L_t^\infty}^2
&\leq
\iint h(w_1) h(w_2) 
|\mathbb{W}(w_1,w_2)|
\, \d w_1 \d w_2
\end{align*}
for $\mathbb W(w_1,w_2)$ defined by \eqref{Wn=proj}.
For a fixed $w_1$, we decompose the space-time domain with respect to $w_2$-variable into translated copies of $\mathbb{V}_{k,l}(w_1)$, shifted vertically in time by
$(0,-2^{-k}\mathit{j})\in \R^n\times \R$ (see also Figure~\ref{f:decomposition}). 
This decomposition yields
\begin{equation}\label{transtang}
    \|h \ast \chi_{\mathbb{V}_{k,l}}\|_{L_x^2L_t^\infty}^2
    \leq
   \mathcal J_{\trans}+ \mathcal J_{\tang},
\end{equation}
where we set
\[
\mathcal J_{\trans}
=\!\!
\sum_{\mathit{j}\ge C 2^{(k-l)\theta}} \int h(w_1) \int \chi_{\mathbb{V}_{k,l}(0,2^{-k}\mathit{j})}(w_2-w_1) h(w_2) |\mathbb{W}(w_1,w_2)|\,\d w_1 \d w_2
\]
and
\begin{align}\label{tangtang}
\mathcal J_{\tang}
=\!\!
\sum_{\mathit{j}\le C 2^{(k-l)\theta}} \int h(w_1) \iint \chi_{\mathbb{V}_{k,l}(0,2^{-k}\mathit{j})}(w_2-w_1) h(w_2) |\mathbb{W}(w_1,w_2)|\,\d w_1 \d w_2
\end{align}
for a constant $C>0$ and parameter $0<\theta<1$ to be chosen later. We refer to these terms as \textit{transversal} and \textit{tangential}, respectively.\footnote{The terminology reflects the geometric nature of the intersection between two cones $\mathbb V_{k,l}(w_1)$ and translation of $\mathbb V_{k,l}(w_2)$ for $w_2 \in \mathbb V_{k,l}(w_1-(0,2^{-k}\mathit{j}))$. For each fixed $w_1$, these intersections exhibit distinct behaviors depending on the size of $\mathit{j}$. For large $\mathit{j}$ ($\ge C2^{(k-l)\theta}$), they intersect more transversely, while for small $\mathit{j}$ ($\le C2^{(k-l)\theta}$), they meet more tangentially. While the transversal case may include some tangential intersections, we maintain this terminology for clarity of presentation.}

\begin{lemma}\label{lem:trta}
For the transversal case, we have
\begin{equation}\label{e:two cones refined2}
|\mathbb{W}(w_1,w_2)|
\lesssim
2^{-\frac{(1+\theta)}{2}k}2^{-\frac{(4-\theta)}{2}l}
(|x_1-x_2|+2^{-k})^{-\frac12},
\quad~ \mathit{j}\ge C2^{(k-l)\theta}
\end{equation}
and for the tangential case, we have
\begin{equation}\label{e:two cones refined1}
|\mathbb{W}(w_1,w_2)|
\lesssim
2^{-\frac k2} 2^{-2l}(|x_1-x_2|+2^{-k})^{-\frac12},
\quad\quad\quad\,\,\, \mathit{j}\le C 2^{(k-l)\theta}.
\end{equation}
\end{lemma}
The detailed computation leading to these bounds can be found in Subsection~\ref{subsec:inter cone}.
Assuming \eqref{e:two cones refined2} and \eqref{e:two cones refined1}, we now proceed to handle the tangential and transversal cases separately.

\medskip

Let us first deal with the transversal case $\mathit{j}\ge C 2^{(k-l)\theta}$, where the cones $\mathbb{V}_{k,l}(w_1)$ and $\mathbb{V}_{k,l}(w_2)$ are sufficiently well separated for any choice of $w_2$, resulting in a fairly small intersection.
Indeed, 
by using \eqref{e:two cones refined2}, bringing the summation inside the integrals so that $\sum_{\mathit{j}} \chi_{\mathbb V_{k,l}(0,2^{-k}\mathit{j})}\lesssim \chi_{\mathbb B^{3}(0,2^{-l})}$, we have
\begin{align}\label{trantran}
\mathcal{J}_\trans
&\lesssim
2^{-\frac{1+\theta}{2}k}2^{-\frac{4-\theta}{2}l}
\iint h(w_1)h(w_2) \chi_{\mathbb B^3(0,2^{-l})}(x_2-x_1)|x_1-x_2|^{-\frac12}\,\d w_1 \d w_2.
\end{align}
Hence by Young's convolution inequality given that $|x_1-x_2|\lesssim 2^{-l}$, we have
\begin{align}\label{trans-bound}
\mathcal{J}_\trans
&\lessapprox
2^{-\frac{(1+\theta)}{2}k}2^{-\frac{(7-\theta)}{2}l}
\|h\|_{L_x^2L_t^1}^2.
\end{align}

\medskip

Next, we estimate $\mathcal{J}_\tang$. As before, we apply \eqref{e:two cones refined1} to the right-hand side of \eqref{tangtang}.
We write the integrand $h(w_1)h(w_2)\chi_{\mathbb{V}_{k,l}(0,2^{-k}\mathit{j})}(w_2-w_1) |x_1-x_2|^{-\frac12}$ by splitting it as $ (h(w_1)h(w_2))^\frac34 \chi_{\mathbb V_{k,l}(0,2^{-k}\mathit{j})}(w_2-w_1)$ and $(h(w_1)h(w_2))^\frac14 |x_1-x_2|^{-\frac12}$.
Applying H\"older's inequality with exponents $p=4/3$ and $p'=4$, note that the second term can be handled similarly as in the case of transversal case (cf. \eqref{trantran}). While the transversal case benefits from the bound $|x_1-x_2|^{-\frac12}$ yielding an additional decay $2^{-\frac 32l}$, the tangential case does not provide such additional decay.

Nevertheless, we can show that \eqref{tangtang} is bounded by
\begin{align*}
2^{\frac k2}2^{2l}\mathcal{J}_\tang
\lessapprox
\sum_{\mathit{j}\leq C2^{(k-l)\theta}}
\Big(\int h(w_1)  \big(h \ast \chi_{\mathbb{V}_{k,l}}\big)(w_1-(0,2^{-k}\mathit{j}))\,\d w_1\Big)^\frac34
\Big(\|h\|_{L_x^2L_t^1}^2\Big)^\frac14.
\end{align*}
By applying H\"older's inequality again, and then using the translation invariance of the Lebesgue measure to eliminate the dependence on $\mathit{j}$ in the integrals, and allowing the summation to be $\sim 2^{(k-l)\theta}$, we get
\begin{align}\label{tang-dominant}
2^{\frac k2}2^{2l}\mathcal{J}_\tang
   &\lessapprox
    2^{(k-l)\theta}
    \Big(\|h\|_{L_x^2L_t^1}\Big)^{\frac24+\frac34}
    \|h*\chi_{\mathbb{V}_{k,l}}\|_{L_x^2L_t^\infty}^\frac34.
\end{align}

In the case where the tangential part is dominant ($\mathcal J_{\trans} \lesssim \mathcal J_{\tang}$), substituting \eqref{tang-dominant} into \eqref{transtang} and rearranging yields the estimates
\[
\|h*\chi_{\mathbb{V}_{k,l}}\|_{L_x^2L_t^\infty}^2
\lessapprox
\left(2^{-(\frac12-\theta)k}2^{-(2+\theta)l}\right)^\frac85
\|h\|_{L_x^2L_t^1}^2.
\]
Conversely, when the transversal part is dominant ($\mathcal J_{\tang} \lesssim \mathcal J_{\trans}$), by \eqref{trans-bound} we have:
\[
\|h*\chi_{\mathbb{V}_{k,l}}\|_{L_x^2L_t^\infty}^2
\lessapprox
2^{-\frac{(1+\theta)}{2}k}2^{-\frac{(7-\theta)}{2}l}
\|h\|_{L_x^2L_t^1}^2.
\]
By optimizing $\theta$ as 
\[
\theta
=
\frac17,
\]
the right-hand sides of the two expressions coincide, yielding the desired bound \eqref{goal2} and hence \eqref{goal1}. \qed

\subsection{The intersection of two cones}\label{subsec:inter cone}
It remains to analyze the intersection of two cones in each tangential and transversal case and establish Lemma~\ref{lem:trta}.
To estimate the intersection between the cones, we slice the space by taking level sets 
\[
\mathbb{L}_m=\{(x,t):|t-c_0(m+m_0-1)2^{-k}|\leq 2^{-k+\frac12}\}
\]
for $m\in\Z$, $c_0\in\R$, and $m_0$ is chosen so that $|\mathbb{L}_1\cap \mathbb{V}_{k,l}(w_1) \cap \mathbb{V}_{k,l}(w_2)|\not = \emptyset$ but $|\mathbb{L}_{0}\cap \mathbb{V}_{k,l}(w_1) \cap \mathbb{V}_{k,l}(w_2)| = \emptyset$. 
For suitable choices of integers $j \in \mathbb{Z}$, 
depending on the relative position of the cones, and after adjusting $c_0$ appropriately, 
one obtains the following two identities:
\begin{equation}
\begin{aligned}\label{r1r2}
|x_1 - x_2| &= r_1(m) + r_2(m) - 2^{-k + \frac{3}{2}} m,\\
|x_1 - x_2| &= r_1(m) - r_2(m) + 2^{-k + \frac{1}{2}} j.
\end{aligned}
\end{equation}	
Solving the system \eqref{r1r2} yields explicit expressions for the radii:
\begin{align*}
r_1(m)
&=
\Big(m-\frac {\mathit{j}}2\Big)2^{-k+\frac12}+|x_1-x_2|,\quad r_2(m)=\Big(m+\frac {\mathit{j}}2 \Big) 2^{-k+\frac12}.
\end{align*}

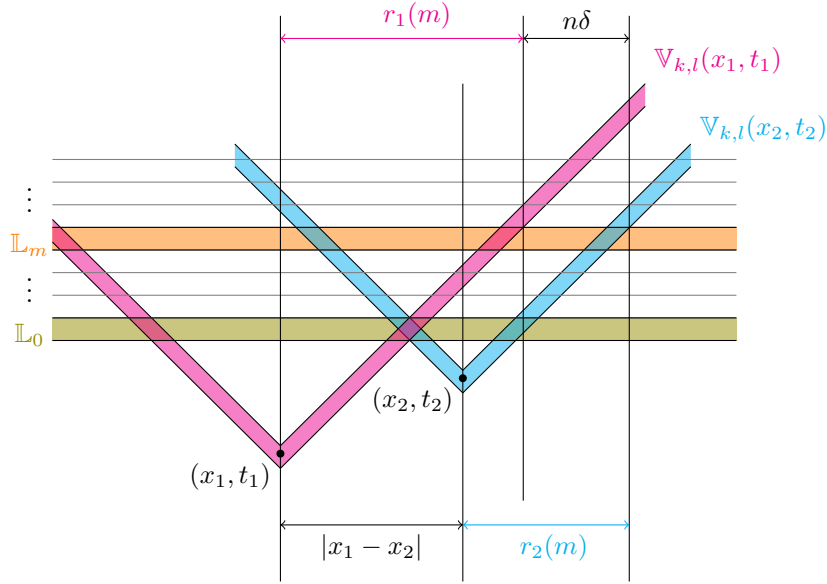
\begin{figure}[htbp]

\begin{tikzpicture}[scale=3]

\fill [olive, opacity=0.5, shift={(0,1/2+0.065)}](-1,0)--(2,0)--(2,0.1)--(-1,0.1);
\fill[orange, opacity=0.5, shift={(0,0.4)}] (-1,1/2+0.065)--(2,1/2+0.065)--(2,1/2+0.165)--(-1,1/2+0.165);

\begin{scope}[shift={(4/5,1/3)}]
    \fill [cyan,opacity=0.5](1,1)--(0,0)--(-1,1)--(-1,1+0.1)--(0,0+0.1)--(1,1+0.1);
    \draw (1,1)--(0,0)--(-1,1);
    \draw[shift={(0,0.1)}] (1,1)--(0,0)--(-1,1);
\end{scope}

\fill [magenta,opacity=0.5](1.6,1.6)--(0,0)--(-1,1)--(-1,1+0.1)--(0,0+0.1)--(1.6,1.6+0.1);
\draw (1.6,1.6)--(0,0)--(-1,1);
\draw[shift={(0,0.1)}] (1.6,1.6)--(0,0)--(-1,1);

\node [magenta,right] at (1.6,1.8) {$\mathbb{V}_{k,l}(x_1,t_1)$}; 
\node [cyan,right] at (1+4/5,1+1/2) {$\mathbb{V}_{k,l}(x_2,t_2)$};

\foreach \y in {1,2,...,8}{
\begin{scope}[gray]
\draw[shift={(0,0.\y65)}] (-1,1/2)--(2,1/2);
\end{scope}
}

\draw[shift={(0,0.065)}] (-1,1/2)--(2,1/2);
\draw[shift={(0,0.165)}] (-1,1/2)--(2,1/2);

\draw[shift={(0,0.465)}] (-1,1/2)--(2,1/2);
\draw[shift={(0,0.565)}] (-1,1/2)--(2,1/2);

\draw (0,2)--(0,-1/2);
\draw (4/5,1.7)--(4/5,-1/2);

\draw [shift={(1.53,0)}](0,2)--(0,-1/2);
\draw [shift={(1.065,0)}](0,2)--(0,-1/7);

\draw [<->,black] (0,-1/4)--(4/5,-1/4);
\draw [<->,cyan] (4/5,-1/4)--(1.53,-1/4);

\draw [<->,magenta] (0,1.9)--(1.065,1.9);
\draw [<->,black] (1.065,1.9)--(1.53,1.9);

\node [shift={(0,0.265)},olive] at (-1.1,1/2) {$\mathbb{L}_0$};
\node [shift={(0,0.965)}] at (-1.1,1/2) {$\vdots$};
\node [shift={(0,2.165)}] at (-1.1,1/2) {$\vdots$};
\node [shift={(0,1.465)}, orange] at (-1.1,1/2) {$\mathbb{L}_m$};

\node [below] at (2/5,-1/4) {$|x_1-x_2|$};
\node [below,cyan] at (1.2,-1/4) {$r_2(m)$};
\node [above,magenta] at (3/5,1.9) {$r_1(m)$};
\node [above] at (1.3,1.9) {$n\delta$};

\fill [black] (0,0+0.065) circle [radius=1/2pt];
\fill [black, shift={(4/5,1/3)}] (0,0+0.065) circle [radius=1/2pt];
\node [below left] at (0,0+0.065) {$(x_1,t_1)$};
\node [below left] at (4/5,1/3+0.065) {$(x_2,t_2)$};

\end{tikzpicture}
\caption{Relationship between $\mathbb{V}_{k,l}(x_1,t_1)$ and $\mathbb{V}_{k,l}(x_2,t_2)$. To estimate the volume of $\mathbb{V}_{k,l}(x_1,t_1)\cap \mathbb{V}_{k,l}(x_2,t_2)$ we further take level sets and apply Wolff's lemma to $\mathbb{O}_\delta$.}
\end{figure}

We first prove \eqref{e:two cones refined1}
Lemma~\ref{l:wolff} yields that 
\[
    |\mathbb{W}(x_1,t_1,x_2,t_2)|
    \lesssim
    \sum_{m} \frac{2^{-\frac{3}{2}k} (r_1(m)+r_2(m))}{|x_1-x_2|^{\frac12}}
    \lesssim
    \frac{2^{-\frac k2} 2^{-2l}}{|x_1-x_2|^\frac12}.
\]
The last inequality follows from the fact that the number of $m$ is $\lesssim 2^{k-l}$ and $r_1(m), r_2(m) \lesssim 2^{-l}$. This immediately yields 
\eqref{e:two cones refined1}.

To show \eqref{e:two cones refined2}, it suffices to consider $t_1\leq t_2$, which implies $r_1(m)\geq r_2(m)$ for all $m$.\footnote{As mentioned earlier, our transversal case still contains tangential cases. When $t_2$ is sufficiently small and distant from $t_1$, two cones are aligned \textit{tangentially} (in this case $\mathbb V_{k,l}(w_1)$ contains the other cone). However, this situation can be handled analogously to $\mathcal J_\tang$.} From the definitions \eqref{d2} and \eqref{r1r2}, we have $\varDelta=(2^{-k+\frac32}m)\cdot(2^{-k+\frac12}\mathit{j})$ and $r_1(m)+r_2(m)=2^{-k+\frac32}m+|x_1-x_2|$. We now set $m_*$ as the smallest value of $m$ satisfying 
\begin{align}\label{mstar}
m\geq 4\frac{|x_1-x_2|}{2^{-k} \mathit{j}}.
\end{align}
With this definition, $\varDelta \ge 2^{-k}  (r_1(m)+r_2(m))$ holds if and only if $m\ge m_*$.
Hence by Lemma~\ref{l:wolff},
\begin{align*}
    |\mathbb{W}(w_1,w_2)|
    &\lesssim
    \sum_{m<m_*} \frac{2^{-\frac32k} ( r_1(m)+r_2(m))}{|x_1-x_2|^{\frac12}}
    +
    \sum_{m\geq m_*} \frac{2^{-2k} ( r_1(m)+r_2(m))^\frac32}{|x_1-x_2|^{\frac12}\varDelta^{\frac12}}\\
    &:=
    \mathcal{W}_1+\mathcal{W}_2.
\end{align*}
Since $r_1(m)+r_2(m)\lesssim |x_1-x_2|$, recalling \eqref{mstar} and using $\mathit{j}\ge C 2^{\theta(k-l)}$, we see
\[
\mathcal{W}_1
\lesssim \frac{|x_1-x_2|}{2^{-k} \mathit{j}}\cdot 2^{-\frac 32k} |x_1-x_2|^{\frac12}
\lesssim |x_1-x_2|^{\frac 32} 2^{-\frac k2} 2^{-\theta(k-l)} \lesssim 2^{-(\frac12+\theta)k} 2^{-(\frac 32-\theta)l}.
\]
For the last inequality, we use $|x_1-x_2|\lesssim 2^{-l}$.

On the other hand, note that $\varDelta = 2^{-2k}m\cdot\mathit{j} \ge C 2^{-2k}  2^{\theta(k-l)} m$. Since the number of $\{m:m \ge m_*\}$ is trivially bounded by $2^{k-l}$, we have
\[
\mathcal{W}_2
\lesssim \sum_{m\ge m_*}m^{-\frac12}\frac{2^{-k}2^{-\frac\theta2 (k-l)}2^{-\frac 32 l}}{|x_1-x_2|^{\frac12}}
\lesssim 
\frac{2^{-\frac{1+\theta}{2}k}2^{-(2-\frac \theta 2)l}}{|x_1-x_2|^{\frac12}}.
\]
Putting together the estimates for $\mathcal{W}_1$ and $\mathcal{W}_2$, we obtain
\[
|\mathbb{W}(x_1,t_1,x_2,t_2)|
\lesssim
    2^{-(\frac12+\theta)k} 2^{-(\frac32-\theta)l} + \frac{2^{-\frac{1+\theta}{2}k} 2^{-(2-\frac{\theta}{2})l}}{|x_1-x_2|^{\frac12}}.
\]
Since the second term on the right-hand side dominates the other, we obtain \eqref{e:two cones refined2}
as desired.

\begin{remark}
    This improvement argument, which relies on the decomposition using the parameter $\theta$, is effective only when $n=2$. In higher dimensions $n\ge3$, the intersection volume $|\mathbb{O}_\delta(x_1,r_1) \cap \mathbb{O}_\delta(x_2,r_2)|$ for transversal case becomes significantly larger than the tangential case, as is shown in Lemma~\ref{l:wolff higher dim}. More precisely, when $\mathit{j}>C2^{\theta(k-l)}$, the transversal contribution dominates the tangential case regardless of the $\theta$-decomposition. This explains why this argument fails to yield better bounds when $n\ge3$.
\end{remark}

\begin{appendix}

\section{Maximal estimates for single datum}\label{appendix-single}
In this section, we review classical maximal estimates
for single wave operators, providing detailed proofs. While these results are well-known (cf. \cite{CHL, HKL}), the proof offers valuable insight into techniques used to prove analogs for orthonormal systems.

Consider initial data in the inhomogeneous Sobolev spaces $H^s$. Since the lower frequency plays no significant role in our analysis, our argument extends naturally to a general family of operators, including Klein--Gordon operators. 
We define $U_t^\alpha f(x) = e^{it\sqrt{\alpha-\Delta}}f(x) = (2\pi)^{-n}\int e^{i(x\cdot \xi + t\sqrt{\alpha+|\xi|^2})} \widehat{f}(\xi)\,\d \xi$, where $U_t^0=U_t$ corresponds to the wave operator. For $n\geq 1$, we have the following results.
\begin{proposition}
Let $s\in[0,\frac n2)$, $r\ge2$, and $\alpha\geq0$. Then
\begin{align}\label{maximal single modified}
\big\|U_t^\alpha f\big\|_{L_x^r(\B^n)L_t^\infty(0,1)}
\lesssim \|f\|_{\dot{H}^s}
\end{align}
holds if 
\[s>\max\Big\{\frac n2-\frac nr,~ \frac{n+1}r-\frac{n-1}{2r}\Big\}.
\]
Furthermore, the regularity orders are sharp up to endpoints.
\end{proposition}
This is typically stated with $\|f\|_{H^s}$ on the right-hand side, but it can be replaced by $\|f\|_{\dot{H}^s}$ due to \eqref{e:lowfrequency single} provided that $s<\frac n2$. 
The model case is, of course, the wave operator when $\alpha=0$ (\cite{Cowling, RV08, CHL}), but a straightforward modification allows us to extend the result to general $\alpha>0$. 
We present the proof of this extension below.

\subsection{Sufficiency}
Let us show the sufficiency part of \eqref{maximal single modified} first. Recall the frequency localization $P_k$ associated with $\varphi_k$ for $k\in\mathbb Z$. Set $\pst = \frac{2(n+1)}{n-1}$. Note that the choice $r=\pst$ is equivalent to $\frac n2-\frac nr=\frac{n+1}{4}-\frac{n-1}{2r}$ whenever $n\geq2$.
We have the trivial bound 
\[
\|U_t^\alpha P_k f\|_{L_x^\infty(\mathbb{B}^n)L_t^\infty(0,1)}
\lesssim
2^{\frac n2 k}
\|P_k f\|_{L^2}, \quad k \in \mathbb Z. 
\]
In particular, for $k<0$ this provides
\begin{equation}\label{e:lowfrequency single}
\sum_{k<0} \|U_t^\alpha f\|_{L_x^\infty(\B^n)L_t^\infty(0,1)}
\lesssim
\|f\|_{\dot{H}^s}
\end{equation}
whenever $s<\frac n2$. By H\"older's inequality,
this trivially implies $\dot{H}^s$--$L_x^pL_t^\infty$ bounds for $U_t^\alpha (\sum_{k<0} P_kf)$.

If we further establish 
\begin{align}
&\|U_t^\alpha P_k f\|_{L_x^2(\mathbb{B}^n)L_t^\infty(0,1)}
\lesssim
2^\frac k2 
\|P_k f\|_{L^2},\label{e:single L2}\\
&\|U_t^\alpha P_k f\|_{L_x^\pst(\mathbb{B}^n)L_t^\infty(0,1)}
\lesssim
2^{(\frac n2-\frac{n}{\pst})k} 
\|P_k f\|_{L^2},\label{e:single Lpst}
\end{align}
for $k\geq0$, a routine interpolation argument yields \eqref{maximal single modified}.

To derive \eqref{e:single L2} and \eqref{e:single Lpst}, we use the Sobolev embedding, which follows directly from the fundamental theorem of calculus, and H\"older's inequality. Indeed, for $r\ge1$,
\begin{align*}
\sup_{0<t<1} |U_t^\alpha P_kf(x)|^r
\le
|U_0^\alpha P_k f(x)|^r
+
2^k\Big(\int_0^1|U_t^\alpha P_k f(x)|^{r}\,\d t\Big)^{r-1} \Big(\int_0^1| U_t^\alpha\tilde{P_k}f(x)|^r\,\d t\Big)
\end{align*}
holds where $\tilde P_k$ is a projection operator satisfying $U_t^\alpha \tilde P_k f =2^{-k} \partial_t(U_t^\alpha P_kf)$. Since $\|\tilde P_kf\|_2\sim \|P_kf\|_2$, it suffices to consider $U_t^\alpha P_kf$. 

When $r=2$, by Plancherel's theorem we have
\[
\|U_t^\alpha P_kf\|_{L^2(\R^n \times (0,1))}
\lesssim
\|P_kf\|_{L^2}.
\]
When $r=\pst$, Strichartz estimate (after appropriate rescaling) implies
\[
\|U_t^\alpha P_k f\|_{L^\pst(\R^{n+1})}
\lesssim
2^{(\frac n2-\frac{n+1}{\pst})k}
\|P_kf\|_{L^2}.
\]
Hence, by the Sobolev embedding those estimates give \eqref{e:single L2} and \eqref{e:single Lpst}.
This finishes the proof.
\qed

\subsection{Necessity}
We now prove the necessity part. 
We first show that $s\ge \frac n2-\frac nr$ is necessary for \eqref{maximal single modified} to hold. Take $\widehat f(\xi)=\chi_{\mathbb B^n(0,N)}(\xi)$. Then for $t=1/N$, it follows that $|U_t^\alpha f(x)|\gtrsim N^n$ provided $|x|\le N^{-1}$ and $N\ge \alpha$. Hence $N^{\frac nr}\le \big\|\sup_t|U_t^\alpha f| \big\|_{L_x^r}\lesssim \|f\|_{H^s}\lesssim N^{s+\frac n2}$ holds for large $N$. Taking $N\rightarrow\infty$, we see that $s\ge \frac n2-\frac nr$ is necessary.

\smallskip
To see $s\geq \frac{n+1}{4}-\frac{n-1}{2r}$ is necessary for \eqref{maximal single modified} to hold, we consider the function
\[
\widehat{f}(\xi) 
=
N^{-\frac{n-1}{4}} \phi_P(\xi)
\]
for $P=\{\xi=(\xi_1,\bar{\xi})\in\R\times\R^{n-1}: N\leq \xi_1\leq 2N,\ |\bar{\xi}|\leq c_0N^\frac12\}$ and a sufficiently small $0<c_0\ll 1$. Clearly, $|P|\simeq N^\frac{n+1}{2}$ and $\|f\|_{L^2}\lesssim N^\frac12$.

Rearranging the phase, we express
\[
|U_t^\alpha f(x) |
=
(2\pi)^{-n}N^{-\frac{n-1}{4}}
\Big|\int e^{i(\bar x\cdot \bar \xi +(x_1+t)\xi_1+t(\sqrt{\alpha^2+|\xi|^2}-\xi_1))}
\phi_P(\xi)\,\d \xi \Big|.
\]
Using the fact that $\sqrt{\alpha^2+|\xi|^2}-\xi_1=\xi_1\sqrt{1+(\alpha^2+|\overline \xi|^2)/\xi_1^2}-\xi_1=(\alpha^2+|\overline \xi|^2)/\xi_1+O(c_0^4)$ provided that $\xi \in P$ and $N\ge\alpha$.
One can verify that the phase is bounded by $\frac 12$ when $(x,t)$ belonging to the set $Q=\{(x_1,\bar x, t)\in \R\times\R^{n-1}\times \R: |x_1+t|\leq \frac{1}{100}N^{-1},\ |\bar{x}|\leq \frac{1}{100} N^{-\frac12},\ |t|\leq \frac{1}{100}\}$ with $|\proj_tQ|\simeq N^{-\frac{n-1}{2}}$. Consequently, we obtain
\begin{align*}
    \|U_t^\alpha f\|_{L_x^rL_t^\infty(Q)}
    \gtrsim
    N^{\frac{n+3}{4}} |Q|^\frac1r 
    \simeq
    N^{\frac{n+3}{4}-\frac{n-1}{2r}}.
\end{align*}

Therefore, it follows 
$N^{\frac{n+3}{4}-\frac{n-1}{2r}} \lesssim N^{s+\frac12}$.
Letting $N\to \infty$, one obtains that 
$s\geq \frac{n+1}{4}-\frac{n-1}{2r}$.
\qed

\subsection{Remarks on endpoint regularity}
For applications such as the pointwise convergence problem, it is crucial to determine the least regularity in terms of $r\in[1,\infty)$, meaning that one can choose the optimal $r$ that minimizes the required regularity. This is usually attained when $r=2$ (sometimes it is the only case), and our case here is no exception. In fact, $\min_{r\geq2} \max\{\frac n2-\frac nr,\frac{n+1}{r}-\frac{n-2}{2r}\} = \frac12$ when $r=2$. In this special case, the estimate
\begin{equation}\label{e:single r=2}
\|U_t^\alpha f\|_{L_x^2(\mathbb{B}^n)L_t^\infty((0,1))}
\lesssim
\|f\|_{H^s}
\end{equation}
is known to hold if and only if $s>\frac 12$, as shown by Walther \cite{Walther}.

To see the failure of \eqref{e:single r=2} when $s=\frac12$, let $\widehat \zeta\in C_0^\infty$ be supported on $[\frac12,2]$ with $\int \widehat \zeta(s)s^{\frac{n-1}2}\,ds=1$. For fixed $0<r<1$, we set $\widehat{\zeta_R}(\xi)=R^{-\frac{n+1}2}\widehat \zeta(R^{-1}|\xi|)$ for $R\ge 1/r$. Using polar coordinates, we obtain
\begin{align*}
U_t^\alpha \zeta_R
=R^{-\frac{n+1}2}\int \int_{\mathbb S^{n-1}} e^{i sx\cdot \theta}\,\d \theta \,\,e^{i t\sqrt{\alpha^2+s^2}}\widehat \zeta(R^{-1}s)s^{n-1}\,\d s.
\end{align*}

The asymptotic behavior of the Bessel functions yields
$\int_{\mathbb S^{n-1}}e^{i s x\cdot \theta} \,\d \theta= \linebreak \sum_{\pm}c_{\pm} ( s |x|)^{-\frac{n-1}2} e^{\pm i |x|s}+O\big(s^{-\frac{n+1}2}|x|^{-\frac{n+1}2})$ for some constants $c_\pm\in \R$. Setting $t=|x|$, we have
\begin{align*}
U_{|x|}^\alpha \zeta_R(x)
=c_{-}R^{-\frac{n+1}2}|x|^{-\frac{n-1}2}\int e^{i |x|(\sqrt{\alpha^2+s^2}-s)}\widehat\zeta\Big(\frac sR\Big)s^{\frac{n-1}2}\,\d s+O(R^{-1}|x|^{-\frac{n+1}2}).
\end{align*}
Here, the integral involving the integrand $c_+(|x|/s)^{-\frac{n-1}2}e^{i|x|(\sqrt{\alpha^2+s^2}+s)}$ can be shown to be absorbed into the error term $O(R^{-1}|x|^{-\frac{n+1}2})$ through repeated integration by parts.

\smallskip

For the wave operator ($\alpha=0$), the above integral is  straightforward to obtain. For general $\alpha>0$, Taylor series expansion gives $\sqrt{\alpha^2+s^2}-s=\alpha^2/(2s)+O(\alpha^4/s^3)$. Thus, the phase function in the integral is smaller than $1/2$ provided $|x|\le1$ and for large enough $R\gg \alpha^2$. Therefore, there exists $|c(x)|\sim1$ such that
\[
U_{|x|}^\alpha \zeta_R(x)= c(x)|x|^{-\frac{n-1}2}+O(R^{-1}|x|^{-\frac{n+1}2}).
\]

Now we take $f_\lambda$ such that $\widehat f_\lambda=\sum_{1\le 2^k \le \lambda}\widehat\zeta_{2^k}$.
Then it follows that
\[
U_{|x|}^\alpha f_\lambda (x) = \sum_{1 \le 2^k \le \lambda}c(x)|x|^{-\frac{n-1}2}+O(1).
\]
For $r \le |x|\le 1$ and $\lambda \gtrsim r^{-2}$, we have $\sup_{0<t\le1}|U_t^\alpha f_\lambda| \gtrsim \log \lambda$. Consequently,
$
\big\| \sup_{0<t<1}|U_t^\alpha f_\lambda| \big\|_{L^2} \gtrsim \log \lambda.
$
Since $\|f_\lambda\|_{H^{1/2}}\sim (\log \lambda)^{\frac 12}$, assuming \eqref{e:single r=2} holds leads to a contradiction.  \qed

\section{Proof of Lemma~\ref{l:wolff} and \ref{l:wolff higher dim}}\label{Appendix:B}
We begin with an elementary but crucial observation for the proof of Lemma~\ref{l:wolff}.
\begin{lemma}\label{l:angle size}
Let $-1<\mu<1$ and $\Theta=\{\theta\in [-\frac \pi2,\frac\pi 2]: |\cos \theta - \mu|\le \delta\}$. 
Then
\[
|\Theta| \lesssim \frac{\delta} {\sqrt{|1-\mu|+\delta}}.
\]
\end{lemma}
This estimate follows directly from the Taylor series expansion. We omit the detailed proof.

\subsection{Proof of Lemma~\ref{l:wolff}}
We now prove \eqref{CC refined}, following an argument in \cite[Lemma 3.1]{Wolff}.
We first divide the cases based on the separation of centers and deal with $|x_1-x_2|\le 4\delta$ and $|x_1-x_2|\ge 4\delta$ separately. 

Suppose first that $|x_1-x_2| \le 4\delta$. Then it is clear that
\[
|\mathbb{O}_\delta(x_1,r_1)\cap \mathbb{O}_\delta(x_2,r_2)|\lesssim \delta r_1+\delta r_2,
\]
and this, together with $\Delta \lesssim \delta(r_1+r_2)$, yields \eqref{CC refined}. 

\smallskip

Now we consider the case $|x_1-x_2|\ge 4\delta$.
Without loss of generality, we may assume that $x_1=(0,0)$, $x_2:=(y,0)\in \R \times \R$ with $y>0$, and $r_1>r_2$. In this case, we observe that
$r_1\sim r_1+r_2$, since it is clear that $r_1\le r_1+r_2 \le 2r_1$.
We define the parameter \[d = |x_1-x_2|+|r_1-r_2|,\] which simplifies to $d=y+(r_1-r_2)$.
We claim that
\begin{align}\label{x sim d}
4^{-1}d \le y \le d.
\end{align}

To establish the lower bound, we may assume that
$
y\ge {d}/2-\delta.
$
Otherwise, from the definition of $d$, we would have $y<d/2-\delta=y/2+(r_1-r_2)/2-\delta$, which rearranges to $y+r_2<r_1-2\delta$. For $x_2=(y,0)$, this leads to
$\mathbb{O}_\delta(0,r_1)\cap \mathbb{O}_\delta(x_2,r_2)=\emptyset$.
Since $d\ge 4\delta$, we have $y\ge d/2-\delta \ge d/4$. The upper bound $y\le d$ is clear. This proves the claim \eqref{x sim d}.

\smallskip

Now we prove \eqref{CC refined} under the assumption \eqref{x sim d}. To simplify the analysis, we rescale all parameters by $r_1^{-1}$ and consider the intersection
$\mathbb O_{\delta/r_1}(0,1)\cap \mathbb O_{\delta/r_1}(x_2/r_1,r_2/r_1)$.
Under this scaling, we may assume that the radius is $r_1=1$.

In the normalized setting, let $z\in \mathbb{O}_\delta(0,1)\cap \mathbb{O}_\delta(x_2,r_2)$. Then $|z-e^{i\theta}|\le \delta$ for some $\theta$.
Since $e^{i\theta}\in \mathbb{O}_\delta(x_2,r_2)$, we have $|\,|e^{i\theta}-(y+i0)|^2-r_2^2|\le \delta$, which leads to
$(\cos \theta-y)^2+\sin^2 \theta-r_2^2=O(\delta)$. Consequently, we obtain 
$
1+y^2-r_2^2-2y\cos \theta=O(\delta)
$
from which it follows that
\begin{align*}
\cos \theta 
&=(1+y^2-r_2^2)/2y +O(\delta/y).
\end{align*}

We set $\mu:=(1+y^2-r_2^2)/(2y)$. 
Recalling the definition \eqref{d2}, we may write
\begin{align}\label{mu approx}
1-\mu=
\frac{r_2^2-(1-y)^2}{2y}
=\frac{(r_2-1+y)(1+r_2-y)}{2y}
=\frac{\varDelta}{2y}.
\end{align}

Next, by applying Lemma~\ref{l:angle size} with $\delta$ replaced by $\delta/y$, and taking into account the thickness $\delta$ of the circle, along with the approximation $y \sim y+\delta$ (cf. \eqref{x sim d}) we see
\[
|\mathbb{O}_\delta(0,1)\cap \mathbb{O}_\delta(x_2,r_2)|
\lesssim
\delta \cdot \frac{\delta}{y}\Big(\frac{y}{\varDelta+\delta}\Big)^{\frac12} \lesssim \frac{\delta^2}{(y+\delta)^{\frac12}(\varDelta+\delta)^{\frac12}}. 
\]

For general $r_1$, we replace $(\delta,y,r_2)$ by $r_1^{-1}(\delta,y,r_2)$ (thus replacing $\varDelta$ with $r_1^{-2}\varDelta$), and then dilate the intersection by $r_1$. This gives
\[
|\mathbb{O}_\delta(x_1,r_1)\cap \mathbb{O}_\delta(x_2,r_2)|
\lesssim
\frac{\delta^{\frac 32} r_1}{(y+\delta)^{\frac12}} \left(\frac{\varDelta+\delta r_1}{\delta r_1}\right)^{-\frac12}.
\]
Since $r_1\sim r_1+r_2$ and $y=|x_1-x_2|$, this completes the proof of Lemma~\ref{l:wolff}.
\qed

\medskip

For later use in the proof of Lemma~\ref{l:wolff higher dim}, we note the following result:
\begin{lemma}
Suppose $r_1\geq r_2$.
The intersection $\mathbb{O}_\delta(x_1,r_1)\cap \mathbb{O}_\delta(x_2,r_2)$ is contained within a disc centered at $x_2+r_1 \sgn(r_1-r_2) \frac{x_1-x_2}{|x_1-x_2|}$ with radius approximately
\begin{align}\label{radius size}
\left(\frac{(r_1+r_2)\delta +\varDelta} {|x_1-x_2|+\delta} \right)^\frac12.
\end{align}
\end{lemma}

Indeed, under the assumption that $r_1=1>r_2$, it follows from \eqref{mu approx} that $\theta$ is contained in an interval of length approximately
$
\big((\varDelta+\delta)/(y+\delta)\big)^{1/2}.
$
Rescaling $(\delta,y)$ by $r_1^{-1}(\delta,y)$ (and thus $\varDelta$ by $r_1^{-2}\varDelta$) followed by dilation by $r_1$ as above, yields \eqref{radius size}.

\medskip

\subsection{Proof of Lemma~\ref{l:wolff higher dim}}

We may assume that the centers $x_1,x_2$ lies in a common $2$-dimensional plane in $\R^n$. For simplicity, we write $x_1=(x_1',0)$ and $x_2=(x_2',0)$ with $x_1',x_2' \in \R^2$. If
$x=(x',x'')\in \R^2\times \R^{n-2}$ are contained in $\mathbb{O}_\delta((x_1',0),r_1)\cap \mathbb{O}_\delta((x_2',0),r_2)$, then
\[
x' \in \mathbb{O}_\delta(x_1',r_1)\cap \mathbb{O}_\delta(x_2',r_2), \quad
x''\in \{z''\in \R^{n-2}: |z''|=\sqrt{1-|x'|^2}+O(\delta)\}.
\]
By \eqref{radius size}, the distance of $x''$ from the $2$-plane $\R^2$ is bounded above by a constant multiple of $\big(\frac{\varDelta+\delta (r_1+r_2)}{|x_1-x_2|+\delta}\big)^{1/2}$.
Recalling \eqref{CC refined}, we obtain
\[
|\mathbb{O}_\delta(x_1,r_1)\cap \mathbb{O}_\delta(x_2,r_2)|
\lesssim
\frac{\delta^{\frac32}(r_1+r_2)}{(|x_1-x_2|+\delta)^{\frac12}}
	\Big(\frac{\delta(r_1+r_2)}{|x_1-x_2|+\delta}\Big)^{\frac12}\Big(\frac{\varDelta+\delta (r_1+r_2)}{|x_1-x_2|+\delta}\Big)^{\frac{n-3}2}
\]
Thus we have \eqref{intersect high} as desired.
\qed

\end{appendix}

\subsection*{Acknowledgement} 
This work was supported by JSPS KAKENHI Grant number JP24K16945 (Kinoshita), the NRF (Republic of Korea) grant RS-2024-00339824, G-LAMP RS-2024-00443714,  JBNU research funds for newly appointed professors in 2024, and NRF2022R1I1A1A01055527 (Ko), and by FCT/Portugal through project UIDB/04459/2020 with DOI identifier 10-54499/UIDP/04459/2020 (Shiraki). The third author thanks David Beltran for enlightening discussions. The authors also express their gratitude to Neal Bez, Sanghyuk Lee and Shohei Nakamura for earlier collaboration and valuable discussions.

\end{document}